\documentclass[11pt]{article}

\usepackage{graphicx}
\usepackage{amssymb}
\usepackage{amsmath}

\newtheorem{thm}{Theorem}
\newtheorem{lem}[thm]{Lemma}

\newtheorem{definition}{Definition}
\newtheorem{cor}{Corollary}

\def\ds{\displaystyle}

\title{A Positive Flux Limited Difference Scheme for Option Pricing 2D Fully Non-linear Parabolic Equation with  Uncertain Correlation}
\author{Miglena N. Koleva, Lubin G. Vulkov\\[0.8ex]
\small{\textit{University of Rousse , 8 Studentska St., 7017 Rousse, Bulgaria}}\\[-0.8ex]
\small {\{\textit{mkoleva,lvalkov}\}\textit{@uni-ruse.bg}}}
\date{}
\begin{document}

\maketitle

\begin{abstract} We consider a two-asset non-linear model of option pricing in an environment where the correlation is not known precisely, as it varies between two known values. First we discuss the non-negativity of the solution of the problem. Next, we construct and analyze a positivity preserving, flux-limited finite difference scheme for the corresponding boundary value problem. Numerical experiments are analyzed.\\[-0.1in]


\end{abstract}

\section{Introduction}

Very important for the valuation of  option pricing models  is the correct specification of the respective model parameters. Some of them are given from the market, or  estimated from historic or forward looking data but others are the result of  calibration to market prices. These techniques leads to more realistic in practice \textit{non-linear}  models with uncertain parameter values, for example volatility, interest rate, dividend or correlation.

Usually this parameters range between upper and lower known bonds and consequently we may consider highest and lowest option value, called \textit{best} and \textit{worst} values. These prices can be interpret as \textit{worst-case pricing} for short and long position respectively.

Well-known one-factor uncertain volatility models are derived by Avellaneda, Levy and Par\'{a}s \cite{ALP}. Following Black-Scholes hedging and no-arbitrage arguments they construct  a  worst/best option pricing model where the value of the  volatility depends on the sign of the second derivative, the  Gamma greek ($\Gamma$).

The same idea  applied to the case of uncertain interest rate or uncertain dividend yield (independent of the asset price) in the case of continuous dividend 
 leads to  non-linear one-asset uncertain parameter models, which gives a consistent way to eliminate the dependence of a price on a parameter and to some extent reduce model dependence \cite{Wil1}.

The same arguments \cite[p.313]{Wil1} can be carried over to  multi-asset models, strongly dependent on the correlation $\rho$ between the stochastic processes of the underlying state variable. The correlation  is difficult to guess or calculate in practice so it can be considered as uncertainty. Following \cite{BS} and \cite{Wil1}, this simple hedging strategy is realized in \cite{Top1} for two-asset option pricing model. To be self-contained we  outline the derivation of the model, presented in \cite{Top1}.

Consider the correlation bounded by $-1\leq \rho_1\leq \rho \leq\rho_2\leq 1$ and define the price movements of two underlying assets $S_1$, $S_2$ (for time $t$, trends (drift rates) $\mu_1$, $\mu_2$, volatilities $\sigma_1$, $\sigma_2$ and increments of standard Wiener's process $dX$)
\begin{eqnarray*}
dS_1=\mu_1S_1dt+\sigma_1S_1dX,\\ dS_2=\mu_2S_2dt+\sigma_2S_2dX,
\end{eqnarray*}
correlated by $E(dX_idX_j)=\rho dt$.

By It\^{o}'s Lemma we express an infinitesimal change in the portfolio ($\Pi$), consisting of a long position in one option and short position in both underlyings. Next, eliminating the risk, just as in the classical argument when deriving the Black-Scholes equation for the option prise $V(S_1,S_2,t)$ we get
$$
d\Pi=\left(\frac{\partial V}{\partial t}+\frac12\sigma_1^2S_1^2\frac{\partial^2V}{\partial S_1^2}+\frac12\sigma_2^2S_2^2\frac{\partial^2V}{\partial S_2^2}+\rho\sigma_1\sigma_2S_1S_2\frac{\partial^2V}{\partial S_1\partial S_2}\right)dt.
$$

In order to derive worst-case scenario model we will be  extremely pessimistic: in every infinitesimal time step we assume that a correlation leads to the smallest growth in the portfolio, i.e.
\begin{equation}\label{Pf0}
\min\limits_{\rho} d\Pi=r\Pi\, dt, \ \ \hbox{where} \ \ r>0 \ \ \hbox{is the interest rate.}
\end{equation}

Taking into account that the portfolio  consists of a long position in one option and short position in both underlying we have
\begin{equation}\label{Pf1}
r\Pi \, dt=r\left(V(S_1,S_2,t)-\frac{\partial V}{\partial S_1}S_1- \frac{\partial V}{\partial S_2}S_2\right)dt
\end{equation}
and
\begin{eqnarray}
\min\limits_{\rho} d \Pi&=&\min\limits_{\rho}\left\{\left(\frac{\partial V}{\partial t}+\frac12\sigma_1^2S_1^2\frac{\partial^2V}{\partial S_1^2}+\frac12\sigma_2^2S_2^2\frac{\partial^2V}{\partial S_2^2}+\rho\sigma_1\sigma_2S_1S_2\frac{\partial^2V}{\partial S_1\partial S_2}\right)dt\right\}\nonumber\\[-0.05in]
\label{Pf2}\\[-0.05in]
&=&  \left\{\begin{array}{ll}
\ds\frac{\partial V}{\partial t}+\frac12\sigma_1^2S_1^2\frac{\partial^2V}{\partial S_1^2}+\frac12\sigma_2^2S_2^2\frac{\partial^2V}{\partial S_2^2}+\rho_1\sigma_1\sigma_2S_1S_2\frac{\partial^2V}{\partial S_1\partial S_2}, & \ds \frac{\partial^2V}{\partial S_1\partial S_2}>0,\\
\ds \frac{\partial V}{\partial t}+\frac12\sigma_1^2S_1^2\frac{\partial^2V}{\partial S_1^2}+\frac12\sigma_2^2S_2^2\frac{\partial^2V}{\partial S_2^2}+\rho_2\sigma_1\sigma_2S_1S_2\frac{\partial^2V}{\partial S_1\partial S_2}, & \ds \frac{\partial^2V}{\partial S_1\partial S_2}<0.
\end{array}\right.\nonumber
\end{eqnarray}

Combining \eqref{Pf1}, \eqref{Pf2} via \eqref{Pf0} and taking into account the dividends (denoted by $D_1$ and  $D_2$) we obtain the worst-case pricing equation
\begin{eqnarray}
\begin{split}\label{Eq}
&\frac{\partial V}{\partial t}+\frac12\sigma_1^2S_1^2\frac{\partial^2V}{\partial S_1^2}+\frac12\sigma_2^2S_2^2\frac{\partial^2V}{\partial S_2^2}+\rho(\Gamma_{cross})\sigma_1\sigma_2S_1S_2\frac{\partial^2V}{\partial S_1\partial S_2}\\&\hspace{0.6in}+(r-D_1)S_1\frac{\partial V}{\partial S_1}+(r-D_2)S_2\frac{\partial V}{\partial S_2}-rV=0, \ \ (S_1,S_2)\in \Omega=\mathbb{R}^+\times \mathbb{R}^+, \ \ 0\leq t<T;
\end{split}\\[0.1in]
\rho(\Gamma_{cross})=\left\{\begin{array}{ll}
\rho_1,&\Gamma_{cross}>0,\\
\rho_2,&\Gamma_{cross}<0.
\end{array}\right.,\ \ \ \  \Gamma_{cross}=\frac{\partial^2V}{\partial S_1\partial S_2},\ \ \ \ -1\leq \rho_1\leq \rho_2\leq 1.\hspace{0.2in}\label{Nonl}&
\end{eqnarray}

In the best-case scenario for an investor with long position, $\rho(\Gamma_{cross})$ is determined by
$$
\rho(\Gamma_{cross})=\left\{\begin{array}{ll}
\rho_1,&\Gamma_{cross}<0,\\
\rho_2,&\Gamma_{cross}>0.
\end{array}\right.
$$

There are many numerical methods for {one-asset} uncertain parameter models available in the literature.  For example, for the  uncertain volatility model (which is identical with Leland model of transaction cost \cite{Wil1}), in \cite{PFV} is developed numerical iteration algorithm. Positivity preserving method is presented in \cite{K}.
 A fully-implicit, monotone discretization method is developed for the solution of  option pricing model with uncertain drift rate in \cite{WWFV}.

 For multi-asset (or two-asset) \textit{linear models}, various numerical methods can be found in the literature, e.g. \cite{CJFC}, where the authors present positivity preserving numerical approach for two-asset linear option pricing stochastic volatility model.

Amid numerous publications, related to the numerical solution of option pricing models, the investigations concerning non-linear multi-asset option pricing models are scarce.
The only work (we managed to find in the literature),  related to the non-linear two-asset option pricing model with uncertain correlation,  is the paper of J. Topper \cite{Top1}. The author implement  the collocation finite element method with cubic Hermite trial functions to solve the worst-case scenario for the considered problem.

In \cite{Ma} a two-asset stochastic correlation model is considered, where the correlation coefficient is a random walk following the square root process. This leads to linear model that is solved by quasi-Monte Carlo method.

In this paper we develop a second-order positivity preserving numerical method for the problem \eqref{Eq},\eqref{Nonl}. We construct implicit-explicit difference scheme,  using different stencils, in dependence of the sign of correlation, for the approximation of $\Gamma_{cross}$  and application of van Leer flux limiter approach for the first derivative discretization. Mild restrictions for space and time mesh step sizes  guarantee the stability and \textit{positivity preserving property} of the numerical solution, i.e. starting with non-negative initial data to obtain a non-negative numerical solution at each time layer.

The rest of the paper is organized as follows. In the next section, we formulate  the differential problem on bounded domain,  after application of the exponential variable change \cite{E,TR}. The non-negativity of the solution is discussed. Combining the monotone techniques in \cite{R,Sam1} with flux limiting, we perform a space discretization of the problem in Section 3. A positive fully-discrete scheme is derived in the next section. Numerical results are discussed in Section 5 and the paper is completed by some conclusions.

\section{The differential problem}

Let now $\overline{\Omega}=\Omega\cup \partial\Omega=[L_W,L_E]\times[L_S,L_N]\subseteq \mathbb{R}^+\times \mathbb{R}^+$. Following the financial modelling in \cite{Top1} we consider the  equation \eqref{Eq}, \eqref{Nonl}, 
associated with the  terminal and boundary conditions \cite{Top0,Top1,Top2,Top3}
\begin{eqnarray}
V(S_1,S_2,T)&=&g_0(S_1,S_2)\geq 0\; \hbox{in}\; \Omega \label{TC}\\
\frac{\partial V(S_1,S_2,t)}{\partial n}&=&g_1(S_1,S_2,t)\geq 0 \; \hbox{on}\;  \partial \Omega_1,
\label{BC1}\\
 V(S_1,S_2,t)&=&g_2(S_1,S_2,t)\geq 0 \; \hbox{on}\;  \partial \Omega_2 \not\equiv \emptyset,\ \ \partial \Omega_1\cup\partial \Omega_2=\partial \Omega.\label{BC2}
\end{eqnarray}

Here $\partial/\partial n$ is the outward derivative to $S_1$ or $S_2$ and $T$ is time to maturity.

Using the logarithmic prices 
\begin{equation}\label{ChV}
x_i=\ln S_i,\ \ i=1,2,\ \ \ \ \tau=T-t,
\end{equation}
we introduce the operators
\begin{eqnarray*}
\mathcal{L}_iu&=&-\frac12\sigma_1^2\frac{\partial^2u}{\partial x_1^2}-\frac12\sigma_2^2\frac{\partial^2u}{\partial x_2^2}-\rho_i\sigma_1\sigma_2\frac{\partial^2u}{\partial x_1\partial x_2}\\&&-(r-D_1-\frac12\sigma_1^2)\frac{\partial u}{\partial x_1}-(r-D_2-\frac12\sigma_2^2)\frac{\partial u}{\partial x_2}+ru, \ \ \ \ i=\{0,1,2\},
\end{eqnarray*}
where we formally set $\rho_0=\rho(\widetilde{\Gamma}'_{cross})$. Then \eqref{Eq}-\eqref{BC2} is transformed to
 the following problem for $u(x_1,x_2,\tau)=V(S_1,S_2,t)$, $(x_1,x_2)\in \overline{\Omega'}=[\ln{L_W},\ln{L_E}]\times[\ln{L_S},\ln{L_N}]\subseteq \mathbb{R}^2$.
\begin{eqnarray}
&& \frac{\partial u}{\partial \tau}+\mathcal{L}_0u=0, \ \ (x_1,x_2,\tau)\in Q_T\equiv \Omega'\times (0,T); \label{EqT}\\
 && \Gamma'_{cross}=e^{-(x_1+x_2)}\frac{\partial^2u}{\partial x_1\partial x_2},\ \ \ \ \widetilde{\Gamma}'_{cross}=\frac{\partial^2u}{\partial x_1\partial x_2},\hspace{1in}\label{NonlT}
 \end{eqnarray}\\[-0.4in]
 \begin{eqnarray}
u(x_1,x_2,0)&=&g'_0(x_1,x_2)\; \hbox{in}\; \Omega' \label{IC}\\
\frac{\partial u(x_1,x_2,\tau)}{\partial n'}&=&{g}'_1(x_1,x_2,\tau) \; \hbox{on}\;  \partial \Omega'_1,\label{BC1T}\\
 u(x_1,x_2,\tau)&=&g'_2(x_1,x_2,\tau) \; \hbox{on}\;  \partial \Omega'_2,\ \ \partial \Omega'_1\cup\partial \Omega'_2=\partial \Omega',\label{BC2T}
  \end{eqnarray}
  where $\partial/\partial n'$ is the outward derivative to $x_1$ or $x_2$, ${g}'_0(x_1,x_2)=g_0(e^{x_1},e^{x_2})$, ${g}'_2(x_1,x_2,\tau)=g_2(e^{x_1},e^{x_2},\tau)$ and
$$
{g}'_1(x_1,x_2,\tau)=\left\{\begin{array}{ll}
e^{x_1}g_1(e^{x_1},e^{x_2},\tau),& x_1=\ln L_W \ \ \hbox{or}\ \ x_1=\ln L_E,\\
e^{x_2}g_1(e^{x_1},e^{x_2},\tau),& x_2=\ln L_S\ \ \hbox{or}\ \ x_2=\ln L_N.
\end{array}\right.
$$

The notation $(\cdot)'$ indicates the transformed by \eqref{ChV} object $(\cdot)$.

Due to the complexity of the presented nonlinear model there are difficulties in obtaining existence and
uniqueness results for problem \eqref{EqT}-\eqref{BC2T}. In this paper we are not concerned with this aspect of the
problem but we shall discuss the minimum principle.

We denote by $C^{m,q} (Q_{T})$ the space of functions defined on $Q_{T}$ that
have continuous derivative with respect to  $x=(x_{1},x_{2})$ up to order $m$
and continuous derivative with respect to $t$ up to order $q$.

Typically, no $C^{2,1}$ solution exists on the hole domain
$Q_{T}$ of equation \eqref{EqT} with discontinuous function $\rho _0$. The
particularity of the equation \eqref{EqT} is that it shows degeneracy, because
it is possible
$\widetilde{\Gamma}'_{cross} = 0$. Thus it is
naturally to assume the existence of a set
$S (x_{1} , x_{2} , \tau) \subset Q_{T} $ on which
$ \widetilde{\Gamma}'_{cross} ( x_{1} , x_{2} , \tau) =0$. This set (it is expected to be a surface) is not given
in advance so that we  have a Stefan-like problem. But
\eqref{EqT} is derived from stochastic finance and therefore specific interface
(internal boundary) conditions are needed.
We assume  $ u \in C^{2,1} (Q_{T} )$ across the phase-change surfaces that is in
accordance with condition $ \widetilde{\Gamma}'_{cross} ( x_{1} , x_{2} , \tau) \vert_{\mathcal{S}} = 0 $. Out of the
interface  $ \mathcal{S} ( x_{1} , x_{2} , \tau) $ we assume even higher regularity,
$ u \in C^{3,1} ( \Omega_{T} \backslash  \mathcal{S})$. By $\partial \Omega_T^p$ we denote the parabolic boundary of $\overline{Q}_T$, i.e. $\partial \Omega_T^p=\{(x_1,x_2,\tau):(x_1,x_2)\in\partial \Omega', 0\leq \tau <T\}$, i.e. the boundary of $Q_T$ minus the interior of the top part of the boundary, $\Omega'\times \{\tau=T\}$. Also, by $Q_T^+$ ($Q_T^-$) we will denote the subset of $Q_T$, where $\widetilde{\Gamma}'_{cross} >0$ ($\widetilde{\Gamma}'_{cross} <0$).

\begin{thm}[Minimum Principle] Suppose that the function
$ u \in C ( {\overline{Q}}_{T} ) \cap C^{2,1} ( Q_{T} ) \cap
  C^{3,1} ( \Omega_{T}  \backslash \mathcal{S} ) $ satisfies
  in $Q_{T} $ the problem \eqref{EqT}-\eqref{BC2T} and $g'_0(x_1,x_2)\geq 0$ in $\Omega'$ and $g'_i(x_1,x_2,\tau)\geq 0$ on $\partial\Omega'_i$, $i=1,2$. Then $u$ can not attain negative local minimum in
  $ { \overline {Q}}_{T} \setminus \partial Q_T^p$ and $u\geq 0$ on $\overline{Q}_T$.
\end{thm}
  \noindent \textbf{Proof.} Suppose that there exists a local minimum point $P_0(x_{1_0},x_{2_0},\tau_0)\in {Q}_T$ with $u(P_0)<0$.

  1. If $0<\tau_0<T$, then $P_0$ belongs to the interior of $Q_T$ and therefore,
  \begin{equation}
  \frac{\partial u}{\partial \tau}(P_0)=\frac{\partial u}{\partial x_1}(P_0)=\frac{\partial u}{\partial x_2}(P_0)=0,\label{e15}
    \end{equation}
  and
   \begin{equation}
 \frac{\partial^2 u}{\partial x_1^2}(P_0)\geq 0,\ \ \frac{\partial^2 u}{\partial x_2^2}(P_0)\geq 0.\label{e16}
    \end{equation}

    1.1. Suppose $P_0\in S$. Then $\frac{\partial^2 u}{\partial x_1 \partial x_2}=0$ and \eqref{e15}, \eqref{e16} lead to
    $$
    \left( \frac{\partial u}{\partial \tau}+\mathcal{L}_0u\right)(P_0)<0,
    $$
    which contradicts to equation \eqref{EqT}.

    1.2. Suppose that $P_0\in Q_T^+$ (similar is the treatment of the case $P_0\in Q_T^-$). Then, in view of \eqref{e15}, \eqref{e16} we have
    \begin{eqnarray}
    \begin{split}
    0&=\left(\frac{\partial u}{\partial \tau}+\mathcal{L}_1u\right)(P_0)=\mathcal{L}_1u(P_0)\\ &=-\frac12\sigma_1^2\frac{\partial^2 u}{\partial x_1^2}(P_0)-\frac12\sigma_2^2\frac{\partial^2 u}{\partial x_2^2}(P_0)-\rho_1\sigma_1\sigma_2\frac{\partial^2 u}{\partial x_1 \partial x_2}(P_0)+ru(P_0).\label{e17}
     \end{split}
      \end{eqnarray}

    Since $P_0$ is not on the boundary of $Q_T$, there is a neighborhood of $(x_{1_0},x_{2_0},t_0)$ within of the domain $Q_T$ where we can use the Taylor expansion:
      \begin{eqnarray*}
      & &u(x_{1_0}+\triangle x_1,x_{2_0}+\triangle x_2,\tau_0)=u(P_0)\\ &&+\frac12\left(\frac{\partial^2 u}{\partial x_1^2}(P_0)(\triangle x_1)^2+2\triangle x_1\triangle x_2\frac{\partial^2 u}{\partial x_1 \partial x_2}(P_0)+\frac{\partial^2 u}{\partial x_2^2}(P_0)(\triangle x_2)^2\right)+O((\triangle x_1)^3+(\triangle x_2)^3).
    \end{eqnarray*}
    Taking into account that $u(x_{1_0}+\triangle x_1,x_{2_0}+\triangle x_2,\tau_0)> u(P_0)$ for all $\triangle x_1$ and $\triangle x_2$ that are small enough, we have
       \begin{equation}\label{e18}
       \frac{\partial^2 u}{\partial x_1^2}(P_0)(\triangle x_1)^2+2\triangle x_1\triangle x_2\frac{\partial^2 u}{\partial x_1 \partial x_2}(P_0)+\frac{\partial^2 u}{\partial x_2^2}(P_0)(\triangle x_2)^2\geq 0.
          \end{equation}
          Since $u(P_0)<0$, from \eqref{e17} follows  that
          \begin{equation}\label{e18a}
       \sigma_1^2\frac{\partial^2 u}{\partial x_1^2}(P_0)+2\rho_1\sigma_1\sigma_2\frac{\partial^2 u}{\partial x_1 \partial x_2}(P_0)+\sigma_2^2\frac{\partial^2 u}{\partial x_2^2}(P_0)< 0.
          \end{equation}

          In order to match the Taylor expansion to get a contradiction, we require the last inequality as
           \begin{equation}\label{e19}
       \left(\frac{\sigma_1}{\sqrt{C}}\right)^2\frac{\partial^2 u}{\partial x_1^2}(P_0)+2\rho_1 \frac{\sigma_1}{\sqrt{C}}\frac{\sigma_2}{\sqrt{C}}\frac{\partial^2 u}{\partial x_1 \partial x_2}(P_0)+ \left(\frac{\sigma_2}{\sqrt{C}}\right)^2\frac{\partial^2 u}{\partial x_2^2}(P_0)< 0,
          \end{equation}
          where $C>0$ is a constant. Next we take
          $$
          \triangle x_1=\frac{\sigma_1}{\sqrt{C}}\ \ \hbox{and} \ \ \triangle x_2=\frac{\sigma_2}{\sqrt{C}}.
          $$
          This contradicts to \eqref{e18a} for sufficiently large $C$.

1.3. Suppose $P_0\in \partial\Omega'_1$ and for concreteness let $P_0(\ln L_W,x_2,\tau)$, i.e. ${x_1}_0=\ln L_W$. Then following similar considerations as in the Hopf's lemma \cite{Ev}, we conclude that $\partial u/\partial n(P_0)>0$, where $n(P_0)$ is the outer normal. But $\partial u/\partial n (P_0)=-\partial u/\partial x (P_0)=g_1(P_0)\leq 0$, so we get contradiction.

          2. Now suppose $\tau_0=T$. Then we will have $\frac{\partial u}{\partial \tau}(P_0)\leq 0$, instead of $\frac{\partial u}{\partial \tau}(P_0) =0$ in \eqref{e15} and we once more deduce the contradiction in the cases 1.1, 1.2 and 1.3. \hfill $\Box$

  \section{Space discretization }
 In the present section we develop the numerical method, combining the idea of A. Samarskii et al. \cite{Sam1} to use different stencils for the approximation of the mixed derivative with the flux limiter approach \cite{GGWC,HV,LV} in two space directions for approximation of the first derivatives.

We define an uniform mesh in space $\overline{\Omega}$
\begin{eqnarray*}
&\overline{\omega}_h=\left\{x=({x_1}_i,{x_2}_j):\; {x_1}_i=L_W+(i-1)h_1,\ \ {x_2}_j=L_S+(j-1)h_2,\right.\\
&\left. i=1,\dots,N_1,\; j=1,\dots,N_2,\ \  h_1={(L_W-L_E)}/{(N_1-1)}, \; h_2={(L_N-L_S)}/{(N_2-1)}\right\}
\end{eqnarray*}
and denote the numerical solution at point $({x_1}_{i},{x_2}_{j},\tau)$ by $u_{i,j}(\tau):=u({x_1}_{i},{x_2}_{j},\tau)$.

\begin{center}
 \begin{figure}
{\includegraphics[width=0.32\textwidth,keepaspectratio=true,
    ]{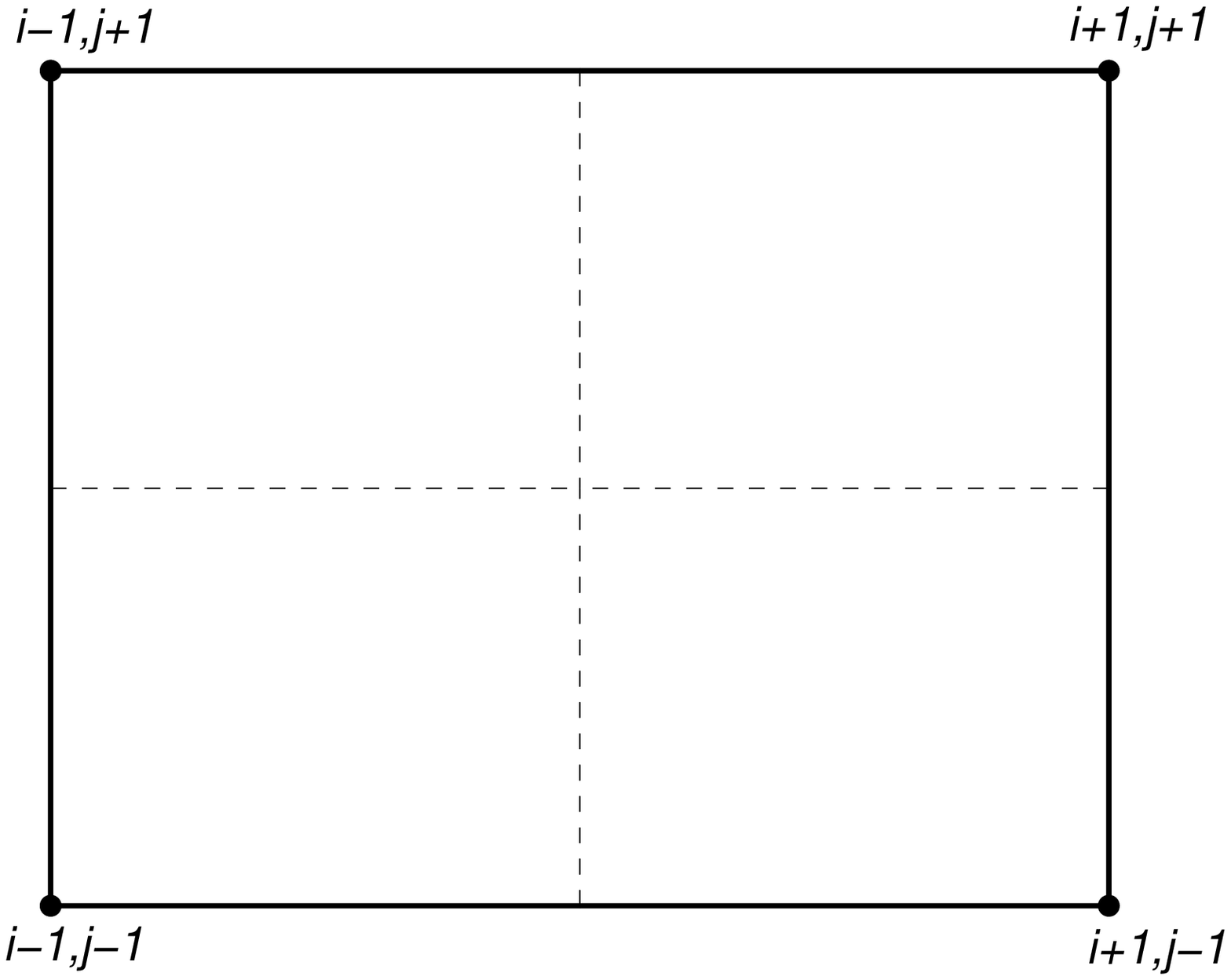}}
\hfill
\includegraphics[width=0.32\textwidth,keepaspectratio=true,
    ]{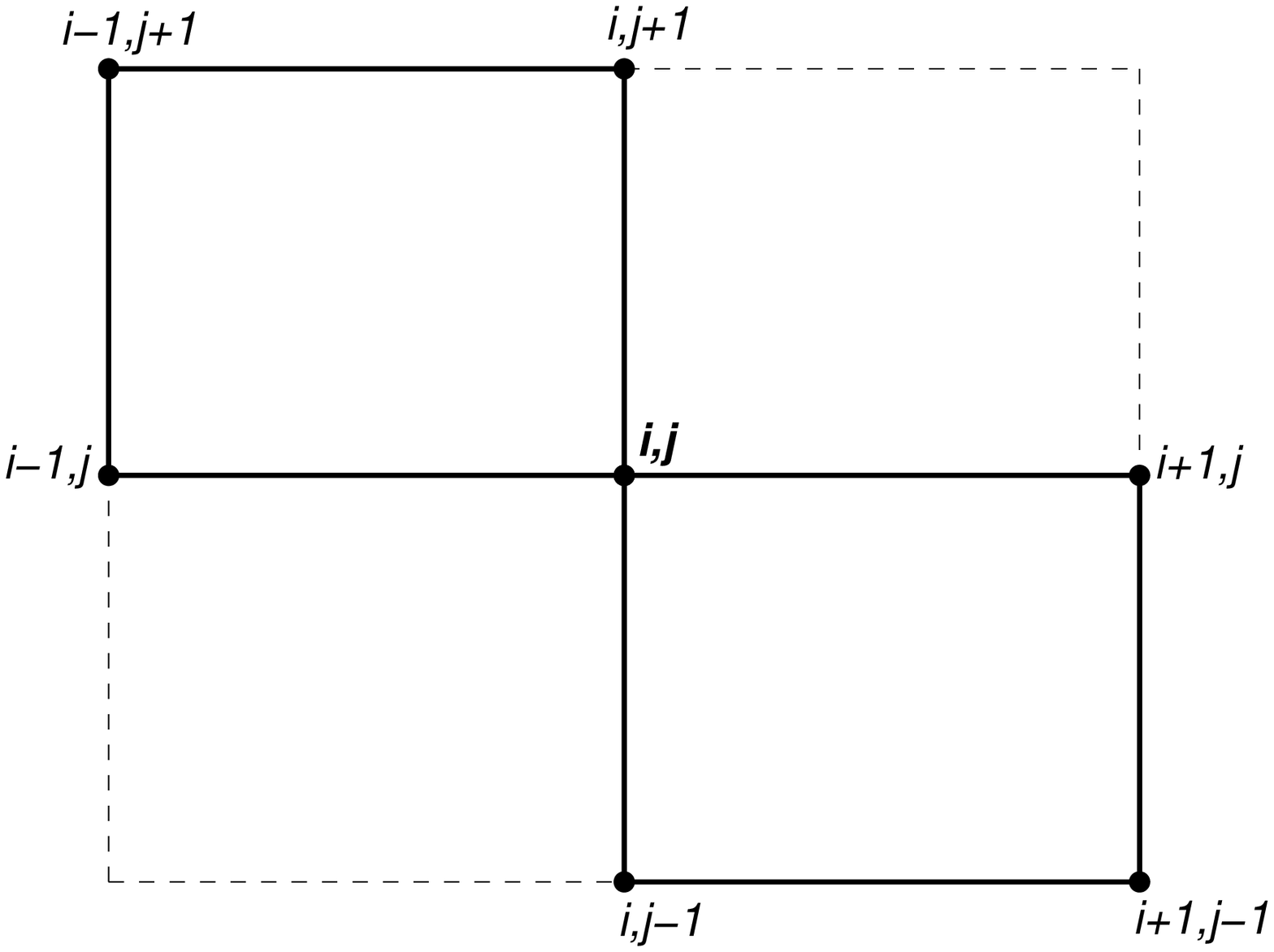}
\hfill
\includegraphics[width=0.32\textwidth,keepaspectratio=true,
    ]{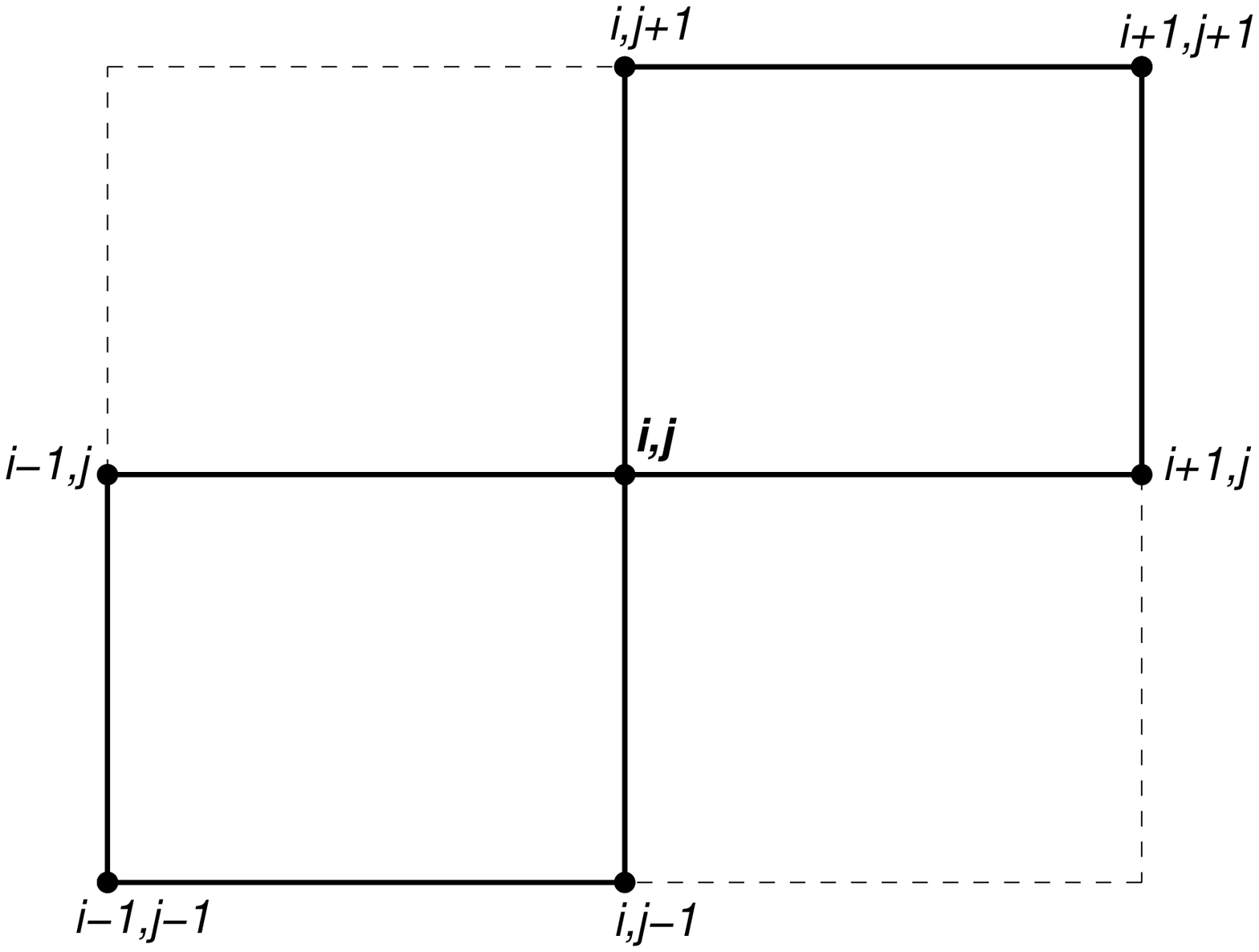}
    \\[-5mm]
\caption{Stencils, corresponding to  ${u_{\mathring{x}_1\mathring{x}_{2_{i,j}}}}$ (\textit{left}), $u_{{x_1x_2}_{i,j}}^-$ (\textit{center}) and $u_{{x_1x_2}_{i,j}}^+$ (\textit{right}) }
\label{f1}
\end{figure}
\end{center}

Further, we use the notations
\begin{eqnarray*}
&&\ds {u_{{\overline{x}_1}_{i,j}}}=\frac{u_{i,j}-u_{i-1,j}}{h_1}, \ \ {u_{{x_1}_{i,j}}}={u_{{\overline{x}_1}_{i+1,j}}}, \ \ {u_{{\overline{x}_2}_{i,j}}}=\frac{u_{i,j}-u_{i,j-1}}{h_2}, \ \ {u_{{x_2}}}_{i,j}={u_{{\overline{x}_2}_{i,j+1}}},  \\
&&\ds u_{\mathring{x}_{s_{i,j}}}=\frac12[{u_{{x_s}_{i,j}}}+{u_{{\overline{x}_s}_{i,j}}}],\ \  {u_{\overline{x}_sx_p}}={(u_{\overline{x}_s})_{x_p}},\ \ {u_{\mathring{x}_s\mathring{x}_p}}=({u_{\mathring{x}_s})_{\mathring{x}_p}}, \ \ s,p \in \mathbb{N},\\
&&\ds  u_{{x_1x_2}_{i,j}}^-=\frac12[u_{{\overline{x}_1}{x_2}_{i,j}}+u_{x_1{{\overline{x}_2}_{i,j}}}], \ \ u_{{x_1x_2}_{i,j}}^+=\frac12[u_{x_1{x_2}_{i,j}}+u_{{\overline{x}_1}\,{\overline{x}_2}_{i,j}}],\ \ \hbox{see Figure \ref{f1}}.
\end{eqnarray*}
We may present an arbitrary function $v$  in the form $v=v^+-v^-$ (and $|v|=v^++v^-$), where $v^+=\max\{0,v\}$ and $v^-=\max\{0,-v\}$. Thus, according to \eqref{Nonl} and \eqref{NonlT} for $\rho'_{i,j}:=\rho(\widetilde{\Gamma}'_{{cross}_{i,j}})$ we have
\begin{equation}\label{rho}
\rho'_{i,j}={\rho'}_{i,j}^+-{\rho'}_{i,j}^-=\left\{\begin{array}{ll}
\rho_1^+-\rho_1^-,& \widetilde{\Gamma}'_{{cross}_{i,j}}> 0,\\
\rho_2^+-\rho_2^-,& \widetilde{\Gamma}'_{{cross}_{i,j}} < 0.\\
\end{array}\right.
\end{equation}

For approximation of the first derivatives in \eqref{EqT} we apply van Leer flux limiter technique \cite{GGWC,HV,LV} in both space directions. Consider the conservative derivatives approximation
\begin{eqnarray}
\begin{split}\label{FD}
&A_s\frac{\partial u}{\partial x_s}=A_s\frac{\partial u}{\partial x_s}\simeq A_s\frac{U_{\mathrm{e}_s+1/2}-U_{\mathrm{e}_s-1/2}}{h_s},\ \    s=\{1,2\},\ \ \hbox{where}\\
&A_s=r-D_s-\frac12 \sigma_s^2, \ \ U_{\mathrm{e}_s\pm q}=\left\{\begin{array}{ll}
u_{i\pm q,j},& s=1,\\
u_{i,j\pm q},& s=2,
\end{array}\right.\ \ q\in \mathbb{R}.
\end{split}
\end{eqnarray}

Using   gradient ratios
\begin{equation}\label{ratio}
\theta_{\mathrm{e}_s+1/2}=\frac{u_{{x_s}_{i,j}}}{u_{{\overline{x}_s}_{i,j}}},  
\end{equation}
we define van Leer flux limiter \cite{GGWC,HV,L74}
\begin{equation}\label{FL}
\Phi(\theta)=\frac{|\theta|+\theta}{1+|\theta|}.
\end{equation}
Observe that
$\Phi(\theta)$ is Lipschitz continuous, continuously differentiable for all $\theta\neq 0$, and
\begin{equation}\label{Pr}
\Phi(\theta) = 0,\ \ \hbox{if} \ \ \theta \leq 0 \ \ \hbox{and} \ \
\Phi(\theta)\leq2 \min\{1, \theta\}.
\end{equation}
Note that at the extreme points of $u$, the slopes $u_{{x_s}_{i,j}}$ and $u_{{\overline{x}_s}_{i,j}}$ have opposite
signs and
$\Phi(\theta_{\mathrm{e}_s+1/2})=0$.

Following \cite{GGWC}  the numerical flux $U_{\mathrm{e}_s+1/2}$ is approximated in a non-linear way
\begin{equation}\label{U1}
U_{\mathrm{e}_s+1/2}=U_{\mathrm{e}_s}+\frac12\Phi(\theta_{\mathrm{e}_s+1/2})(U_{\mathrm{e}_s}-U_{\mathrm{e}_s-1}).
\end{equation}
 Reflecting the indices that appear in $u_{i,j}$ about $i + 1/2$ or $j + 1/2$  yields \cite{GGWC}
\begin{equation}\label{U2}
U_{\mathrm{e}_s+1/2}=U_{\mathrm{e}_s+1}+\frac12\Phi(\theta_{\mathrm{e}_s+3/2}^{-1})(U_{\mathrm{e}_s+1}-U_{\mathrm{e}_s+2}).
\end{equation}
Similarly, the flux $U_{\mathrm{e}_s-1/2}$, corresponding to \eqref{U1} and \eqref{U2} is defined by shifting the index $s$ (i.e. $i$ or $j$).

Using the symmetry property of the flux limiter $\Phi(\theta)=\theta \Phi(\theta^{-1})$ \cite{KT} and \eqref{ratio}, we approximate $A_s\frac{\partial u}{\partial x_s}$ at point $({x_1}_{i},{x_2}_{j},\tau)$, applying \eqref{U1} and \eqref{U2} in dependence of the sign of $A_s=A_s^+-A_s^-$:
\begin{eqnarray}
\begin{split}\label{Flux}
&\hspace{1in} A_s\frac{\partial u}{\partial x_s}\simeq A_s^+\Lambda_s^+u_{x_s}-A_s^-\Lambda_s^-u_{\overline{x}_s},\ \ s=\{1,2\},\\
&\Lambda_s^+=1+\frac12\Phi(\theta_{\mathrm{e}_s+1/2}^{-1})-\frac12\Phi(\theta_{\mathrm{e}_s+3/2}),\ \ \Lambda_s^-=1+\frac12\Phi(\theta_{\mathrm{e}_s+1/2})-\frac12\Phi(\theta_{\mathrm{e}_s-1/2}^{-1}),
\end{split}
\end{eqnarray}
where $0 \leq \Lambda_s^-\leq 2$ and $0 \leq \Lambda_s^+\leq 2$ in view of \eqref{FL}, \eqref{Pr}.

 We implement the idea of \cite{R} so that we use different stencils for the approximation of the second mixed derivative and by \eqref{Flux}, we obtain the following discretization for \eqref{EqT}  at point $({x_1}_{i},{x_2}_{j},\tau)$, $2<i<N_1-1$, $2<j<N_2-1$:
\begin{eqnarray}
\begin{split}\label{EqSD}
&\frac{\partial u}{\partial \tau}-\frac12\sigma_1^2u_{\overline{x}_1x_1}-\frac12\sigma_2^2u_{\overline{x}_2x_2}-\sigma_1\sigma_2(\rho'^+u_{{x_1x_2}}^+-\rho'^-u_{{x_1x_2}}^-)\\
&\hspace{1.3in}-A_1^+\Lambda_1^+u_{x_1}+A_1^-\Lambda_1^-u_{\overline{x}_1}-A_2^+\Lambda_2^+u_{x_2}+A_2^-\Lambda_2^-u_{\overline{x}_2}+ru=0,
\end{split}
\end{eqnarray}
where $\widetilde{\Gamma}'_{{cross}_{i,j}}\simeq {u_{\mathring{x}_s\mathring{x}}}_{p_{i,j}} $ and $\rho'_{i,j}=\rho'({u_{\mathring{x}_s\mathring{x}}}_{p_{i,j}} )$.

For computing the gradient ratio in grid points for $i=\{2,N_1-1\}$ or  $j=\{2,N_2-1\}$  we need  the values of $u_{i,j}$ at the outer grid nodes $({x_1}_{0},{x_2}_{j},\tau)$, $({x_1}_{N_1+1},{x_2}_{j},\tau)$, $({x_1}_{i},{x_2}_{0},\tau)$ and $({x_1}_{i},{x_2}_{N_2+1},\tau)$ for $1<i<N_1$, $1<j<N_2$. Then  the second-order
extrapolation formulas \cite{Sam} will be used
\begin{eqnarray*}
u_{0,j}=3u_{1,j}-3u_{2,j}+u_{3,j}, \ \ u_{N_1+1,j}=3u_{N_1,j}-3u_{N_1-1,j}+u_{N_1-2,j},\\ u_{i,0}=3u_{i,1}-3u_{i,2}+u_{i,3}, \ \ u_{i,N_2+1}=3u_{i,N_2}-3u_{i,N_2-1}+u_{i,N_2-2}.
\end{eqnarray*}

It is trivial to incorporate Dirichlet boundary conditions \eqref{BC2T} on $\partial\Omega'_2$  in the numerical scheme. Thus, only for illustration, we consider the case $\partial\Omega'_1 \equiv \partial\Omega'$, $\partial\Omega'_2 \equiv \emptyset$ and impose \eqref{BC1T} on the whole boundary.

\textit{West boundary $\partial \Omega'_W$: $i=1$, $1<j<N_2$.} From \eqref{BC1T} we have
\begin{equation}\label{WB}
-u_{\mathring{x}_{1_{1,j}}}={g}'_{1_{1,j}}(\tau),\ \  \hbox{and therefore}\ \ u_{0,j}=2h_1{g}'_{1_{1,j}}(\tau)+u_{2,j}, \ \ j=2,\dots,N_2-1.
\end{equation}
Applying \eqref{EqSD} for $i=1$, $1<j<N_2$, where the term $-\Lambda_1^-u_{{\overline{x}_1}_{1,j}}$ is replaced by ${g}'_{1_{1,j}}(\tau)$ and $u_{0,j}$, $u_{0,j\pm 1}$  are eliminated from \eqref{WB}, we get
\begin{eqnarray}
\begin{split}\label{EqSD_W}
&\frac{\partial u}{\partial \tau}-\frac{\sigma_1^2}{h_1}u_{x_1}-\frac12\sigma_2^2u_{\overline{x}_2x_2}-\frac{\sigma_1\sigma_2}{2}|\rho'|(u_{x_1x_2}-u_{x_1{\overline{x}_2}})
-A_1^+\Lambda_1^+u_{{x}_1}-A_2^+\Lambda_2^+u_{x_2}\\
&\hspace{0.2in}+A_2^-\Lambda_2^-u_{\overline{x}_2}+ru=A_1^-{g}'_{1}+\frac{\sigma_1^2}{h_1}g'_1-{\sigma_1}{\sigma_2}({\rho'}^+ g'_{1_{{\overline{x}_2}}}-{\rho'}^-g'_{1_{{x_2}}}),\ \ \  \ \rho'_{1,j}=\rho'({-g'_1}_{{{\mathring{x}_{2_{1,j}}}}}).
\end{split}
\end{eqnarray}

\textit{North boundary $\partial \Omega'_N$: $1<i<N_1$, $j=N_2$.} Now \eqref{BC1T} is replaced by
\begin{equation}\label{NB}
u_{\mathring{x}_{2_{i,N_2}}}={g}'_{1_{i,N_2}}(\tau)\ \   \Rightarrow \ \ u_{i,N_2+1}=2h_2{g}'_{1_{i,N_2}}(\tau)+u_{i,N_2-1}, \ \ i=2,\dots,N_1-1.
\end{equation}
As before, from \eqref{EqSD} at point $(x_{1_{i}},x_{2_{N_2}},\tau)$, replacing $\Lambda_2^+u_{x_{2_{i,N_2}}}$ by ${g}'_{1_{i,N_2}}(\tau)$ we obtain
\begin{eqnarray}
\begin{split}\label{EqSD_N}
&\frac{\partial u}{\partial \tau}-\frac12\sigma_2^2u_{\overline{x}_1x_1}+\frac{\sigma_2^2}{h_2}u_{\overline{x}_2}-\frac{\sigma_1\sigma_2}{2}|\rho'|(u_{\overline{x}_1\overline{x}_2}-u_{x_1{\overline{x}_2}})
-A_1^+\Lambda_1^+u_{{x}_1}+A_1^-\Lambda_1^-u_{\overline{x}_1}\\ &\hspace{0.2in}+A_2^-\Lambda_2^-u_{\overline{x}_2}+ru
=A_2^+{g}'_{1}+\frac{\sigma_2^2}{h_2}g'_1+{\sigma_1}{\sigma_2}({\rho'}^+ g'_{1_{{{x}_1}}}-{\rho'}^-g'_{1_{{\overline{x}_1}}}),\ \ \  \ \rho'_{i,N_2}=\rho'({g'_1}_{{{\mathring{x}_{1_{i,N_2}}}}}).
\end{split}
\end{eqnarray}

\textit{East boundary $\partial \Omega'_E$: $i=N_1$, $1<j<N_2$.} Similarly,  \eqref{BC1T} is discretizied by
\begin{equation}\label{EB}
u_{\mathring{x}_{1_{N_1,j}}}={g}'_{1_{N_1,j}}(\tau)\ \  \hbox{and} \ \ u_{N_1+1,j}=2h_1{g}'_{1_{N_1,j}}(\tau)+u_{N_1-1,j}, \ \ j=2,\dots,N_2-1.
\end{equation}
Thus, from \eqref{EqSD} written at grid node $(x_{1_{N_1}},x_{2_{j}},\tau)$, we get the approximation at east boundary
\begin{eqnarray}
\begin{split}\label{EqSD_E}
&\frac{\partial u}{\partial \tau}+\frac{\sigma_1^2}{h_1}u_{\overline{x}_1}-\frac12\sigma_2^2u_{\overline{x}_2x_2}-\frac{\sigma_1\sigma_2}{2}|\rho'|(u_{\overline{x}_1\overline{x}_2}-u_{\overline{x}_1{{x}_2}})
+A_1^-\Lambda_1^-u_{{\overline{x}}_1}-A_2^+\Lambda_2^+u_{x_2}\\
&\hspace{0.2in}+A_2^-\Lambda_2^-u_{\overline{x}_2}+ru=A_1^+{g}'_{1}+\frac{\sigma_1^2}{h_1}g'_1+{\sigma_1}{\sigma_2}({\rho'}^+ g'_{1_{{{x}_2}}}-{\rho'}^-g'_{1_{{\overline{x}_2}}}),\ \ \  \ \rho'_{N_2,j}=\rho'({g'_1}_{{{\mathring{x}_{2_{N_2,j}}}}}).
\end{split}
\end{eqnarray}

\textit{South boundary $\partial \Omega'_S$: $1<i<N_1$, $j=1$.} Now the corresponding discrete boundary condition in \eqref{BC1T} is
\begin{equation}\label{SB}
u_{\mathring{x}_{2_{i,1}}}={g}'_{1_{i,1}}(\tau)\ \   \Rightarrow \ \ u_{i,0}=2h_2{g}'_{1_{i,1}}(\tau)+u_{i,2}, \ \ i=2,\dots,N_1-1.
\end{equation}
The discretization, corresponding to the south boundary is:
\begin{eqnarray}
\begin{split}\label{EqSD_S}
&\frac{\partial u}{\partial \tau}-\frac12\sigma_1^2u_{\overline{x}_1x_1}-\frac{\sigma_2^2}{h_2}u_{{x}_2}-\frac{\sigma_1\sigma_2}{2}|\rho'|(u_{{x}_1{x}_2}-u_{\overline{x}_1{{x}_2}})
-A_1^+\Lambda_1^+u_{{x}_1}+A_1^-\Lambda_1^-u_{\overline{x}_1}\\ &\hspace{0.2in}-A_2^+\Lambda_2^+u_{{x}_2}+ru
=A_2^-{g}'_{1}+\frac{\sigma_2^2}{h_2}g'_1-{\sigma_1}{\sigma_2}({\rho'}^+ g'_{1_{{{\overline{x}}_1}}}-{\rho'}^-g'_{1_{{{x}_1}}}),\ \ \  \ \rho'_{i,1}=\rho'(-{g'_1}_{{{\mathring{x}_{1_{i,1}}}}}).
\end{split}
\end{eqnarray}

\textit{North-West corner node:} $i=1$, $j=N_2$. Following the same technique as before, we eliminate  artificial grid nodes arise in  \eqref{EqSD} (written at point $i=1$, $j=N_2$), using boundary conditions \eqref{WB} for $j=N_2$ and \eqref{NB} for $i=1$ and replace $A_1^-\Lambda_1^-u_{\overline{x}_1}$ by $A_1^- g'_1$ and $A_2^+\Lambda_2^+u_{{x}_2}$ by $A_2^+ g'_1$. More different is the treatment of the  term $u_{0,N_2+1}$:
$$
u_{0,N_2+1}=\left\{\begin{array}{ll}
u_{2,N_2-1}+2h_2g'_{1_{2,N_2}}+2h_1g'_{1_{1,N_2+1}},&\hbox{applying first}\ \ \eqref{WB}, \ \ \hbox{then}\ \ \eqref{NB},\\
u_{2,N_2-1}+2h_2g'_{1_{0,N_2}}+2h_1g'_{1_{1,N_2-1}},&\hbox{applying first}\ \ \eqref{NB}, \ \ \hbox{then}\ \ \eqref{WB}.
\end{array}\right.
$$
Averaging the above quantities we obtain
\begin{eqnarray*}
u_{0,N_2+1}&=&u_{2,N_2-1}+h_2g'_{1_{2,N_2}}+h_1g'_{1_{1,N_2+1}}+h_2g'_{1_{0,N_2}}+h_1g'_{1_{1,N_2-1}}\\
&&=u_{1,N_2-1}+2h_2g'_{1_{2,N_2}}+2h_1g'_{1_{1,N_2-1}}+2h_1h_2({g'_1}_{{{\mathring{x}_{2}}}}-{g'_1}_{{{\mathring{x}_{1}}}})_{1,N_2}.
\end{eqnarray*}
To compute $\rho'(u_{\mathring{x}_1\mathring{x}_2})$ at grid node  $i=1$, $j=N_2$ we proceed similarly:
$$
u_{\mathring{x}_1\mathring{x}_2}=\left\{\begin{array}{rl}\ds
g'_{1_{\mathring{x}_1}},& \hbox{applying}\ \ \eqref{WB},\\
\ds -g'_{1_{{\mathring{x}_2}}},& \hbox{applying}\ \ \eqref{NB},
\end{array}\right.\ \ u_{\mathring{x}_1\mathring{x}_2}\simeq 0.5(g'_{1_{\mathring{x}_1}}-g'_{1_{{\mathring{x}_2}}})\ \ \Rightarrow \ \ \rho'(u_{\mathring{x}_1\mathring{x}_2})\simeq \rho'(g'_{1_{\mathring{x}_1}}-g'_{1_{{\mathring{x}_2}}}),
$$
as we need only the sign of $u_{\mathring{x}_1\mathring{x}_2}$.

Consequently, the approximation at North-West corner node  is
\begin{eqnarray}
\begin{split}\label{EqSD_WN}
\frac{\partial u}{\partial \tau}-\frac{\sigma_1^2}{h_1}u_{{x}_1}+\frac{\sigma_2^2}{h_2}u_{{\overline{x}}_2}+\sigma_1\sigma_2|\rho'|u_{{x}_1{\overline{x}}_2}
-A_1^+\Lambda_1^+u_{{x}_1}+A_2^-\Lambda_2^-u_{{\overline{x}}_2}+ru\\ =(A_1^-+A_2^+){g}'_{1}
+\left(\frac{\sigma_1^2}{h_1}+\frac{\sigma_2^2}{h_2}\right)g'_1+{\sigma_1}{\sigma_2}{\rho'^+}(g'_{1_{{{{x}}_1}}}-g'_{1_{{{\overline{x}}_2}}})+{\sigma_1}{\sigma_2}{\rho'^-}G_{NW},\ \ \hbox{where} \\ G_{NW}=g'_{1_{{{{x}}_1}}}-g'_{1_{{{\overline{x}}_2}}}+{g'_1}_{{{\mathring{x}_{2}}}}-{g'_1}_{{{\mathring{x}_{1}}}}\ \ \hbox{and} \ \  \rho'_{1,N_2}=\rho'[({g'_1}_{{{\mathring{x}_{1}}}}-{g'_1}_{{{\mathring{x}_{2}}}})_{1,N_2}].
\end{split}
\end{eqnarray}

\textit{North-East corner node:} $i=N_1$, $j=N_2$. From \eqref{EqSD}, \eqref{NB} and \eqref{EB} at point $i=N_1$, $j=N_2$ we get
\begin{eqnarray}
\begin{split}\label{EqSD_NE}
\frac{\partial u}{\partial \tau}+\frac{\sigma_1^2}{h_1}u_{{\overline{x}}_1}+\frac{\sigma_2^2}{h_2}u_{{\overline{x}}_2}-\sigma_1\sigma_2|\rho'|u_{{\overline{x}}_1{\overline{x}}_2}
+A_1^-\Lambda_1^-u_{{\overline{x}}_1}+A_2^-\Lambda_2^-u_{{\overline{x}}_2}+ru\\ =(A_1^++A_2^+){g}'_{1}
+\left(\frac{\sigma_1^2}{h_1}+\frac{\sigma_2^2}{h_2}\right)g'_1-{\sigma_1}{\sigma_2}{\rho'^+}G_{NE}-{\sigma_1}{\sigma_2}{\rho'^-}(g'_{1_{{{{\overline{x}}}_1}}}+g'_{1_{{{\overline{x}}_2}}}),\ \ \hbox{where} \\ G_{NE}=g'_{1_{{{{\overline{x}}}_1}}}+g'_{1_{{{\overline{x}}_2}}}-{g'_1}_{{{\mathring{x}_{2}}}}-{g'_1}_{{{\mathring{x}_{1}}}}\ \ \hbox{and} \ \  \rho'_{1,N_2}=\rho'[({g'_1}_{{{\mathring{x}_{1}}}}+{g'_1}_{{{\mathring{x}_{2}}}})_{1,N_2}].
\end{split}
\end{eqnarray}

\textit{South-East corner node:} $i=N_1$, $j=1$. Again, from \eqref{EqSD}, \eqref{EB} and \eqref{SB} at point $i=N_1$, $j=1$ we have
\begin{eqnarray}
\begin{split}\label{EqSD_SE}
\frac{\partial u}{\partial \tau}+\frac{\sigma_1^2}{h_1}u_{{\overline{x}}_1}-\frac{\sigma_2^2}{h_2}u_{{{x}}_2}+\sigma_1\sigma_2|\rho'|u_{{\overline{x}}_1{{x}}_2}
+A_1^-\Lambda_1^-u_{{\overline{x}}_1}-A_2^+\Lambda_2^+u_{{{x}}_2}+ru\\ =(A_1^++A_2^-){g}'_{1}
+\left(\frac{\sigma_1^2}{h_1}+\frac{\sigma_2^2}{h_2}\right)g'_1-{\sigma_1}{\sigma_2}{\rho'^+}(g'_{1_{{{\overline{x}}_1}}}-g'_{1_{{{{{x}}}_2}}})-{\sigma_1}{\sigma_2}{\rho'^-}G_{SE},\ \ \hbox{where} \\ G_{SE}=g'_{1_{{{{\overline{x}}}_1}}}-g'_{1_{{{{x}}_2}}}-{g'_1}_{{{\mathring{x}_{1}}}}+{g'_1}_{{{\mathring{x}_{2}}}}\ \ \hbox{and} \ \  \rho'_{N_1,1}=\rho'[({g'_1}_{{{\mathring{x}_{2}}}}-{g'_1}_{{{\mathring{x}_{1}}}})_{N_1,1}].
\end{split}
\end{eqnarray}

\textit{South-West corner node:} $i=j=1$. As before, from \eqref{EqSD}, \eqref{WB} and \eqref{SB} at point $i=1$, $j=1$ we obtain
\begin{eqnarray}
\begin{split}\label{EqSD_SW}
\frac{\partial u}{\partial \tau}-\frac{\sigma_1^2}{h_1}u_{{{x}}_1}-\frac{\sigma_2^2}{h_2}u_{{{x}}_2}-\sigma_1\sigma_2|\rho'|u_{{{x}}_1{{x}}_2}
-A_1^+\Lambda_1^+u_{{{x}}_1}-A_2^+\Lambda_2^+u_{{{x}}_2}+ru\\ =(A_1^-+A_2^-){g}'_{1}
+\left(\frac{\sigma_1^2}{h_1}+\frac{\sigma_2^2}{h_2}\right)g'_1+{\sigma_1}{\sigma_2}{\rho'^+}(g'_{1_{{{{x}}_1}}}+g'_{1_{{{{{x}}}_2}}})+{\sigma_1}{\sigma_2}{\rho'^-}G_{SW},\ \ \hbox{where} \\ G_{SW}=g'_{1_{{{{{x}}}_1}}}+g'_{1_{{{{x}}_2}}}-{g'_1}_{{{\mathring{x}_{1}}}}-{g'_1}_{{{\mathring{x}_{2}}}}\ \ \hbox{and} \ \  \rho'_{1,1}=\rho'[(-{g'_1}_{{{\mathring{x}_{1}}}}-{g'_1}_{{{\mathring{x}_{2}}}})_{1,1}].
\end{split}
\end{eqnarray}


Now, we are going to investigate conditions, which guarantee the positivity preserving property of the semi-discrete problem.  Further we need the following well known results.

Consider the initial value problem (IVP) for the ODE system
\begin{equation}
u'(\tau)=g(\tau,u(\tau)), \ \ \tau\geq \tau_0,\ \  u(\tau_0)=u^0, \ \ \tau_0\in \mathbb{R},\ \  u^0\in R^p, \ \ g: \mathbb{R}\times \mathbb{R}^p\rightarrow \mathbb{R}^p\label{IVP}
\end{equation}


\begin{definition}[\cite{GGWC}, Positive ODE system, positive semi-discretization] The ODE in \eqref{IVP} and the IVP \eqref{IVP} are said to be positive if $g$ is continuous and \eqref{IVP} has a unique solution for all $\tau_0$ and for all $u^0$, and $u(\tau)\geq 0$ holds for all $ \tau\geq\tau_0$ whenever $u^0\geq 0$. A semi-discretization of a given PDE (with non-negative solution) is called positive if it leads to a positive ODE system.
\end{definition}

\begin{lem}[\cite{H}] \label{L1} Let $g$ is continuous and \eqref{IVP} has a unique solution for all $\tau_0$ and for all $u_0$. The initial value problem \eqref{IVP} is positive if and only if
\[
v_i=0, \ \ v_j\geq 0 \ \ for \; all\ \  j\neq i\ \  \Rightarrow \ \ g_i(\tau,v)\geq 0,
\]
holds for all $\tau$ and any vector $v \in \mathbb{R}^p$ and all $i=1,\dots, p$.
\end{lem}
As a consequence of Lemma \ref{L1} is
\begin{cor}\textnormal{\textbf{(\cite[p. 34]{H1})}} \label{C1} A linear system $u'(\tau)=Au(\tau)$, $A=\{a_{i,j}\}$ is positive iff $a_{i,j}\geq 0$ for all $i\neq j$.
\end{cor}
Guided by this results, we can apply (just as in \cite{GGWC})  the statement of Lemma \ref{L1} and Corollary \ref{C1} for the numerical discretization of  of \eqref{EqT}-\eqref{BC2T}, written  in the form
\begin{eqnarray}
\frac{d u}{d\tau}=C_{i+1,j}u_{i+1,j}+C_{i-1,j}u_{i-1,j}+C_{i,j+1}u_{i,j+1}+C_{i,j-1}u_{i,j-1}+C_{i+1,j-1}u_{i+1,j-1}\nonumber\\+C_{i-1,j-1}u_{i-1,j-1}
+C_{i-1,j+1}u_{i-1,j+1}+C_{i+1,j+1}u_{i+1,j+1}\label{PDS}\\-C_{i,j}u_{i,j}+g(\tau),\ \ \ \  i=1,\dots,N_1,\ \  j=1,\dots,N_2.\nonumber
\end{eqnarray}

\begin{lem} \label{L2} The ODE system, defined by \eqref{PDS} is positive, if all coefficients $C_{\Sigma_{i,j}}=\{C_{i\pm1,j}, C_{i,j\pm1}$, $C_{i\pm1,j\pm1}\}$ are non-negative and $g(\tau)\geq 0$.
\end{lem}
\noindent \textbf{Proof.}  The results follows from Lemma \ref{L1}.\hfill $\Box$
\begin{thm}\label{Th1} The numerical discretization \eqref{EqSD}, combined with Dirichlet boundary conditions (on $\partial \Omega'_2$) and approximations  \eqref{EqSD_W}, \eqref{EqSD_N}, \eqref{EqSD_E}, \eqref{EqSD_S} and \eqref{EqSD_WN}, \eqref{EqSD_NE}, \eqref{EqSD_SE}, \eqref{EqSD_SW} of the Neumann boundary conditions, depending on the boundary $\partial \Omega'_1$, is positive, if
\begin{eqnarray}
\begin{split}\label{RH}
\frac{\sigma_1}{\sigma_2}\max\limits_{{1+b_W\leq i\leq N_1-b_E}\atop{1+b_S\leq j\leq N_2-b_N}}|\rho'| \leq\frac{h_1}{h_2}\leq\frac{\sigma_1}{\sigma_2\max\limits_{{1+b_W\leq i\leq N_1-b_E}\atop{1+b_S\leq j\leq N_2-b_N}}|\rho'|},\ \ \hbox{where}\\
b_Q=\left\{\begin{array}{ll}
1,& \partial \Omega'_Q \subseteq \partial \Omega'_2,\\
0, & \hbox{elsewhere}
\end{array}\right., \ \ Q=\{W,E,N,S\}.
\end{split}
\end{eqnarray}
\end{thm}

\noindent \textbf{Proof.} First we consider the discretization \eqref{EqSD} at inner points: $2<i<N_1-1$, $2<j<N_2-1$. Taking into account that $|\rho'_{i,j}|=\rho'^+_{i,j}+\rho'^-_{i,j}$, the coefficients, corresponding to \eqref{PDS} are
\begin{eqnarray*}
&\ds C_{i\pm1,j}=\frac{\sigma_1^2}{2h_1^2}-\frac{\sigma_1\sigma_2|\rho'_{i,j}|}{2h_1h_2}+\frac{A_{1}^\pm\Lambda_{1_{i,j}}^\pm}{h_1},\\
&\ds C_{i,j\pm1}=\frac{\sigma_2^2}{2h_2^2}- \frac{\sigma_1\sigma_2|\rho'_{i,j}|}{2h_1h_2}+\frac{A_{2}^\pm\Lambda_{2_{i,j}}^\pm}{h_2},\\
&\ds C_{i-1,j+1}=C_{i+1,j-1}=\frac{\sigma_1\sigma_2\rho'^-}{2h_1h_2}, \ \ \ \  C_{i-1,j-1}=C_{i+1,j+1}=\frac{\sigma_1\sigma_2\rho'^+}{2h_1h_2}, \ \ g\equiv 0.
\end{eqnarray*}
To ensure the condition  of Lemma \ref{L2} we require
\begin{equation}\label{Es1}
\frac{\sigma_1}{\sigma_2}\max\limits_{{1<i<N_1}\atop{1<j<N_2}}|\rho'| \leq\frac{h_1}{h_2}\leq\frac{\sigma_1}{\sigma_2^2\max\limits_{{1<i<N_1}\atop{1<j<N_2}}|\rho'|}.
\end{equation}
For equation, corresponding to Neumann condition imposed on the  East boundary ($i=N_1$, $1<j<N_2$) from \eqref{EqSD_E} we have
\begin{eqnarray*}
&\ds C_{N_1-1,j}=\frac{\sigma_1^2}{h_1^2}-\frac{\sigma_1\sigma_2|\rho'_{N_1,j}|}{h_1h_2}+\frac{A_{1}^-\Lambda_{1_{N_1,j}}^-}{h_1},\\
&\ds C_{N_1,j\pm1}=\frac{\sigma_2^2}{2h_2^2}- \frac{\sigma_1\sigma_2|\rho'_{N_1,j}|}{2h_1h_2}+\frac{A_{2}^\pm\Lambda_{2_{N_1,j}}^\pm}{h_2},\\
&\ds C_{N_1-1,j \pm 1}=\frac{\sigma_1\sigma_2|\rho'_{N_1,j}|}{2h_1h_2}, \ \ \ \   g_{N_1,j}=A_1^+{g}'_{1_{N_1,j}}+\frac{\sigma_1^2}{h_1}g'_{1_{N_1,j}}+{\sigma_1}{\sigma_2}({\rho'}^+ g'_{1_{{{x}_2}}}-{\rho'}^-g'_{1_{{\overline{x}_2}}})_{N_1,j}.
\end{eqnarray*}
It is easy to verify that $C_{\Sigma_{N_1,j}}\geq 0$ and $g_{i,N_2}\geq 0$ if
\begin{equation}\label{Es2}
\frac{\sigma_1}{\sigma_2}\max\limits_{{1<j<N_2}}|\rho'_{N_1,j}| \leq\frac{h_1}{h_2}\leq\frac{\sigma_1}{\sigma_2\max\limits_{{1<j<N_2}}|\rho'_{N_1,j}|}.
\end{equation}

Similarly, from \eqref{EqSD_W}, \eqref{EqSD_N}, \eqref{EqSD_S}, corresponding to Neumann boundary condition on $\partial \Omega'_{\{W,S,N\}}$ respectively, to guarantee that $C_{\Sigma_{\partial \Omega'_{\{W,S,N\}}}}\geq 0$ and $g_{\partial \Omega'_{\{W,S,N\}}}\geq 0$, we obtain the estimates
\begin{eqnarray}
\begin{split}\label{Es3}
&\hspace{1.5in}\frac{\sigma_1}{\sigma_2}\max\limits_{{1<j<N_2}}|\rho'_{1,j}| \leq\frac{h_1}{h_2}\leq\frac{\sigma_1}{\sigma_2\max\limits_{{1<j<N_2}}|\rho'_{1,j}|}, \\& \frac{\sigma_1}{\sigma_2}\max\limits_{{1<i<N_1}}|\rho'_{i,1}| \leq\frac{h_1}{h_2}\leq\frac{\sigma_1}{\sigma_2\max\limits_{{1<i<N_1}}|\rho'_{i,1}|},\ \ \frac{\sigma_1}{\sigma_2}\max\limits_{{1<i<N_1}}|\rho'_{i,N_2}| \leq\frac{h_1}{h_2}\leq\frac{\sigma_1}{\sigma_2\max\limits_{{1<i<N_1}}|\rho'_{i,N_2}|}
\end{split}
\end{eqnarray}

Similar estimate is obtained from the discretizations at the corner node, where the  two Neumann boundaries intersects. For example, let  $\{\partial \Omega'_N, \partial \Omega'_E \}\subseteq \partial \Omega'_1$, then from \eqref{EqSD_NE} for all elements of $C_{\Sigma_{N_1,N_2}}$ and $g_{N_1,N_2}$ we have
\begin{eqnarray*}
&\ds C_{N_1-1,N_2}=\frac{\sigma_1^2}{h_1^2}-\frac{\sigma_1\sigma_2|\rho'_{N_1,N_2}|}{h_1h_2}+\frac{A_{1}^-\Lambda_{1_{N_1,N_2}}^-}{h_1},\\
&\ds C_{N_1,N_2-1}=\frac{\sigma_2^2}{h_2^2}- \frac{\sigma_1\sigma_2|\rho'_{N_1,N_2}|}{h_1h_2}+\frac{A_{2}^-\Lambda_{2_{N_1,N_2}}^-}{h_2},\ \
 C_{N_1-1,N_2- 1}=\frac{\sigma_1\sigma_2|\rho'_{N_1,N_2}|}{h_1h_2}, \\ &\ds  g_{N_1,N_2}=\left(A_1^++A_2^++\frac{\sigma_1^2}{h_1}+\frac{\sigma_2^2}{h_2}\right){g}'_{1_{N_1,N_2}}
-{\sigma_1}{\sigma_2}[{\rho'^+}G_{NE}+{\rho'^-}(g'_{1_{{{{\overline{x}}}_1}}}+g'_{1_{{{\overline{x}}_2}}})]_{N_1,N_2}.
\end{eqnarray*}
The requirement  $C_{\Sigma_{N_1,N_2}}\geq 0$ and $g_{N_1,N_2}\geq 0$ leads to the estimate
\begin{eqnarray}
\begin{split}\label{Es4}
\frac{\sigma_1}{\sigma_2}|\rho'_{N_1,N_2}| \leq\frac{h_1}{h_2}\leq\frac{\sigma_1}{\sigma_2|\rho'_{N_1,N_2}|}.
\end{split}
\end{eqnarray}

Similarly, from \eqref{EqSD_WN}, \eqref{EqSD_SE}, \eqref{EqSD_SW} we get
\begin{equation}\label{Es5}
\frac{\sigma_1}{\sigma_2}|\rho'_{1,N_2}| \leq\frac{h_1}{h_2}\leq\frac{\sigma_1}{\sigma_2|\rho'_{1,N_2}|}, \ \ \frac{\sigma_1}{\sigma_2}|\rho'_{N_1,1}| \leq\frac{h_1}{h_2}\leq\frac{\sigma_1}{\sigma_2|\rho'_{N_1,1}|},\ \
\frac{\sigma_1}{\sigma_2}|\rho'_{1,1}| \leq\frac{h_1}{h_2}\leq\frac{\sigma_1}{\sigma_2|\rho'_{1,1}|}.
\end{equation}
Collecting all  results \eqref{Es1}-\eqref{Es5}, we obtain \eqref{RH}. \hfill $\Box$

\section{Full discretization}

In this section we develop an implicit-explicit second-order numerical algorithm which preserves the positivity property of the solution. A semi-implicit and implicit  method are used for the diffusion (the non-linear term is computed at the old time level) and reaction terms respectively while the convection term is approximated explicitly.

The grid points over the time interval $[0,T]$ are defined by $\tau_{n}=\tau_{n-1}+\triangle\tau$, $n=1,2\dots$, $\tau_0=0$. Approximations of $u(x_i,y_j,\tau_n)$ is denoted by $u_{i,j}^n$, but further for simplicity, we use the notations ${\widehat{u}}_{i,j}:=u^n_{i,j}$ and ${u}_{i,j}:=u^{n-1}_{i,j}$, $\widehat{u}_t:=(\widehat{u}-u)/\triangle\tau$

The full discretization of \eqref{EqSD} is
\begin{eqnarray}
\begin{split}\label{EqSD_FD}
&\widehat{u}_t-\frac12\sigma_1^2\widehat{u}_{\overline{x}_1x_1}-\frac12\sigma_2^2\widehat{u}_{\overline{x}_2x_2}-\sigma_1\sigma_2(\rho'^+\widehat{u}_{{x_1x_2}}^+-\rho'^-\widehat{u}_{{x_1x_2}}^-)+r\widehat{u}
=A_1^+\Lambda_1^+u_{x_1}-A_1^-\Lambda_1^-u_{\overline{x}_1}\\ &+A_2^+\Lambda_2^+u_{x_2}-A_2^-\Lambda_2^-u_{\overline{x}_2},\ \ i=2,\dots,N_1-1, \ \  j=2,\dots,N_2-1.
\end{split}
\end{eqnarray}

For non-homogeneous  Neumann boundaries \eqref{BC1T} (if any) we obtain from \eqref{EqSD_W}, \eqref{EqSD_N}, \eqref{EqSD_E},\eqref{EqSD_S}, the following discretization
\begin{eqnarray}
\begin{split}\label{EqSD_W_FD}
&\widehat{u}_t-\frac{\sigma_1^2}{h_1}\widehat{u}_{x_1}-\frac12\sigma_2^2\widehat{u}_{\overline{x}_2x_2}-\frac{\sigma_1\sigma_2}{2}|\rho'|(\widehat{u}_{x_1x_2}-\widehat{u}_{x_1{\overline{x}_2}})+r\widehat{u}=
A_1^+\Lambda_1^+u_{{x}_1}+A_2^+\Lambda_2^+u_{x_2}\\
&-A_2^-\Lambda_2^-u_{\overline{x}_2}+A_1^-{\widehat{g}}'_{1}+\frac{\sigma_1^2}{h_1}\widehat{g}'_1-{\sigma_1}{\sigma_2}({\rho'}^+ \widehat{g}'_{1_{{\overline{x}_2}}}-{\rho'}^-\widehat{g}'_{1_{{x_2}}}),\ \ \ i=1,\ \ j=2,\dots,N_2-1.
\end{split}\\[0.1in]
\begin{split}\label{EqSD_N_FD}
&\widehat{u}_t-\frac12\sigma_2^2\widehat{u}_{\overline{x}_1x_1}+\frac{\sigma_2^2}{h_2}\widehat{u}_{\overline{x}_2}-\frac{\sigma_1\sigma_2}{2}|\rho'|(\widehat{u}_{\overline{x}_1\overline{x}_2}-\widehat{u}_{x_1{\overline{x}_2}})+r\widehat{u}
=A_1^+\Lambda_1^+u_{{x}_1}-A_1^-\Lambda_1^-u_{\overline{x}_1}\\&-A_2^-\Lambda_2^-u_{\overline{x}_2}+
A_2^+{\widehat{g}}'_{1}+\frac{\sigma_2^2}{h_2}\widehat{g}'_1+{\sigma_1}{\sigma_2}({\rho'}^+ \widehat{g}'_{1_{{{x}_1}}}-{\rho'}^-\widehat{g}'_{1_{{\overline{x}_1}}}),\ \  i=2,\dots,N_1-1, \ \  j=N_2.
\end{split}
\end{eqnarray}
\begin{eqnarray}
\begin{split}\label{EqSD_E_FD}
&\widehat{u}_t+\frac{\sigma_1^2}{h_1}\widehat{u}_{\overline{x}_1}-\frac12\sigma_2^2\widehat{u}_{\overline{x}_2x_2}-\frac{\sigma_1\sigma_2}{2}|\rho'|(\widehat{u}_{\overline{x}_1\overline{x}_2}-\widehat{u}_{\overline{x}_1{{x}_2}})+r\widehat{u}=
-A_1^-\Lambda_1^-u_{{\overline{x}}_1}+A_2^+\Lambda_2^+u_{x_2}\\
&-A_2^-\Lambda_2^-u_{\overline{x}_2}+A_1^+{\widehat{g}}'_{1}+\frac{\sigma_1^2}{h_1}\widehat{g}'_1+{\sigma_1}{\sigma_2}({\rho'}^+ \widehat{g}'_{1_{{{x}_2}}}-{\rho'}^-\widehat{g}'_{1_{{\overline{x}_2}}}),\ \ i=N_1, \ \ j=2,\dots,N_2-1.
\end{split}\\[0.1in]
\begin{split}\label{EqSD_S_FD}
&\hspace{-0.08in}\widehat{u}_t-\frac12\sigma_1^2\widehat{u}_{\overline{x}_1x_1}-\frac{\sigma_2^2}{h_2}\widehat{u}_{{x}_2}-\frac{\sigma_1\sigma_2}{2}|\rho'|(\widehat{u}_{{x}_1{x}_2}-\widehat{u}_{\overline{x}_1{{x}_2}})+r\widehat{u}
=A_1^+\Lambda_1^+u_{{x}_1}-A_1^-\Lambda_1^-u_{\overline{x}_1}\\ &+A_2^+\Lambda_2^+u_{{x}_2}
+A_2^-{\widehat{g}}'_{1}+\frac{\sigma_2^2}{h_2}\widehat{g}'_1-{\sigma_1}{\sigma_2}({\rho'}^+ \widehat{g}'_{1_{{{\overline{x}}_1}}}-{\rho'}^-\widehat{g}'_{1_{{{x}_1}}}),\ \ i=2,\dots,N_1-1,\ \ j=1.
\end{split}
\end{eqnarray}

Finally, for the corner nodes, where  the two Neumann boundaries intersects, from \eqref{EqSD_WN}, \eqref{EqSD_NE}, \eqref{EqSD_SE}, \eqref{EqSD_SW} we have
\begin{eqnarray}
\begin{split}\label{EqSD_WN_FD}
&\hspace{-0.65in}\widehat{u}_t-\frac{\sigma_1^2}{h_1}\widehat{u}_{{x}_1}+\frac{\sigma_2^2}{h_2}\widehat{u}_{{\overline{x}}_2}+\sigma_1\sigma_2|\rho'|\widehat{u}_{{x}_1{\overline{x}}_2}+r\widehat{u}
=A_1^+\Lambda_1^+u_{{x}_1}-A_2^-\Lambda_2^-u_{{\overline{x}}_2}+(A_1^-+A_2^+){\widehat{g}}'_{1}\\&
+\left(\frac{\sigma_1^2}{h_1}+\frac{\sigma_2^2}{h_2}\right)\widehat{g}'_1+{\sigma_1}{\sigma_2}{\rho'^+}(\widehat{g}'_{1_{{{{x}}_1}}}-\widehat{g}'_{1_{{{\overline{x}}_2}}})+{\sigma_1}{\sigma_2}{\rho'^-}\widehat{G}_{NW},\ \ i=1,\ \ j=N_2.
\end{split}\\[0.1in]
\begin{split}\label{EqSD_NE_FD}
&\ \ \; \widehat{u}_t+\frac{\sigma_1^2}{h_1}\widehat{u}_{{\overline{x}}_1}+\frac{\sigma_2^2}{h_2}\widehat{u}_{{\overline{x}}_2}-\sigma_1\sigma_2|\rho'|\widehat{u}_{{\overline{x}}_1{\overline{x}}_2}+r\widehat{u}
=-A_1^-\Lambda_1^-u_{{\overline{x}}_1}-A_2^-\Lambda_2^-u_{{\overline{x}}_2}+(A_1^++A_2^+){\widehat{g}}'_{1}\\&\hspace{0.7in}
+\left(\frac{\sigma_1^2}{h_1}+\frac{\sigma_2^2}{h_2}\right)\widehat{g}'_1-{\sigma_1}{\sigma_2}{\rho'^+}\widehat{G}_{NE}-{\sigma_1}{\sigma_2}{\rho'^-}(\widehat{g}'_{1_{{{{\overline{x}}}_1}}}+\widehat{g}'_{1_{{{\overline{x}}_2}}}), \ \ i=N_1,\ \ j=N_2.
\end{split}\\[0.1in]
\begin{split}\label{EqSD_SE_FD}
&\widehat{u}_t+\frac{\sigma_1^2}{h_1}\widehat{u}_{{\overline{x}}_1}-\frac{\sigma_2^2}{h_2}\widehat{u}_{{{x}}_2}+\sigma_1\sigma_2|\rho'|\widehat{u}_{{\overline{x}}_1{{x}}_2}+r\widehat{u}
=-A_1^-\Lambda_1^-u_{{\overline{x}}_1}+A_2^+\Lambda_2^+u_{{{x}}_2}+(A_1^++A_2^-){\widehat{g}}'_{1}\\ &\hspace{0.7in}
+\left(\frac{\sigma_1^2}{h_1}+\frac{\sigma_2^2}{h_2}\right)\widehat{g}'_1-{\sigma_1}{\sigma_2}{\rho'^+}(\widehat{g}'_{1_{{{\overline{x}}_1}}}-\widehat{g}'_{1_{{{{{x}}}_2}}})-{\sigma_1}{\sigma_2}{\rho'^-}\widehat{G}_{SE},\ \ i=N_1,\ \ j=1.
\end{split}
\end{eqnarray}
\begin{eqnarray}
\begin{split}\label{EqSD_SW_FD}
&\hspace{-0.75in}\widehat{u}_t-\frac{\sigma_1^2}{h_1}\widehat{u}_{{{x}}_1}-\frac{\sigma_2^2}{h_2}\widehat{u}_{{{x}}_2}-\sigma_1\sigma_2|\rho'|\widehat{u}_{{{x}}_1{{x}}_2}+r\widehat{u}
=A_1^+\Lambda_1^+u_{{{x}}_1}+A_2^+\Lambda_2^+u_{{{x}}_2} +(A_1^-+A_2^-){\widehat{g}}'_{1}\\&
+\left(\frac{\sigma_1^2}{h_1}+\frac{\sigma_2^2}{h_2}\right)\widehat{g}'_1+{\sigma_1}{\sigma_2}{\rho'^+}(\widehat{g}'_{1_{{{{x}}_1}}}+\widehat{g}'_{1_{{{{{x}}}_2}}})+{\sigma_1}{\sigma_2}{\rho'^-}\widehat{G}_{SW},\ \ i=1,\ \ j=1.
\end{split}
\end{eqnarray}

Next, we  discuss positivity preserving property and stability of the numerical solution.

The system \eqref{EqSD_FD}, associated with Dirichlet boundary conditions and the discretization \eqref{EqSD_W_FD}-\eqref{EqSD_SW_FD}, in the case of Neumann boundary can be written in the following compact form
\begin{eqnarray}
&\ds -C_{i+1,j}\widehat{u}_{i+1,j}-C_{i-1,j}\widehat{u}_{i-1,j}-C_{i,j+1}\widehat{u}_{i,j+1}-C_{i,j-1}\widehat{u}_{i,j-1}-C_{i+1,j-1}\widehat{u}_{i+1,j-1}\nonumber\\[-0.08in]
& \label{CF}\\[-0.08in]
& \ds -C_{i-1,j-1}\widehat{u}_{i-1,j-1}
-C_{i-1,j+1}\widehat{u}_{i-1,j+1}-C_{i+1,j+1}\widehat{u}_{i+1,j+1}+C_{i,j}\widehat{u}_{i,j}=f_{i,j},\nonumber
\end{eqnarray}
 for $i=1,\dots,N_1$, $j=1,\dots,N_2$ and equivalent matrix form
\begin{eqnarray*}
& \ds\mathcal{M}\widehat{U}=\mathcal{F},\ \ \hbox{where}\\& \ds U=[\underbrace{u_{1,1},u_{2,1},\dots,u_{N_1,1}}_{j=1},\dots,\underbrace{u_{1,j},u_{2,j},\dots,u_{N_1,j}}_{2\leq j\leq N_2-1},\dots,\underbrace{u_{1,N_2},u_{2,N_2},\dots,u_{N_1,N_2}}_{j=N_2}]^T,
\end{eqnarray*}
where $\mathcal{M}=\{m_{k,p}\}$ is a square $N_1N_2\times N_1N_2$ matrix  and  $\mathcal{F}=\{f_k\}$, $k=i+(j-1)N_1$ is a column-vectors with $N_1N_2$ known from the previous time level entries.

Following Corollary 3.20 \cite[p.91]{V}, if $\mathcal{M}$ is diagonal dominant matrix with $m_{k,p}\leq 0$ for all $k\neq p$ and $m_{k,k}>0$ for all $1\leq k\leq N_1N_2$, then $M^{-1}>0$. Thus, if  $\mathcal{F}\geq0$, we can conclude that $\widehat{U}\geq 0$. On this base we can prove the following statement
\begin{thm} If $g_s\geq 0$, $s=0,1,2$, \eqref{RH} holds and
\begin{equation}\label{RT}
 \triangle \tau\leq \frac{h_1h_2}{2(|A_{1}|h_2+|A_{2}|h_1)},
\end{equation}
  then the numerical solution of the problem \eqref{EqT}-\eqref{BC2T} (respectively \eqref{Eq}-\eqref{BC2}), obtained by \eqref{EqSD_FD}, associated with Dirichlet boundary conditions and discretization \eqref{EqSD_W_FD}-\eqref{EqSD_SW_FD} (depending on $\partial\Omega$) is non-negative.
\end{thm}

\noindent\textbf{Proof.} We apply induction method: the statement holds for $\tau_0=0$, assume that  it holds at time $\tau_{n-1}$ and prove that this statement holds at time $\tau_{n}$. Thus, via to the time integration, the corresponding assertion  holds at each time level. Let $u^{n-1}\geq 0$.

First, using the compact form \eqref{CF} of the presented numerical scheme,  we  show that $\mathcal{M}^{-1}>0$, which means that matrix $\mathcal{M}$ posses the above mentioned property, i.e. for all $i=1,\dots,N_1$ and $j=1,\dots,N_2$:

\textsf{P1}. $\mathcal{M}$ is diagonally dominant, which is equivalent to $|C_{i,j}|\geq \hspace{-0.3in}\sum\limits_{C_{i+s_1,j+s_2} \in C_{\sum_{i,j}}}\hspace{-0.3in}|C_{i+s_1,j+s_2}|$;\\[-0.05in]

\textsf{P2}. $m_{k,p}\leq 0$ for all $k\neq p$, equivalently to $C_{{i+s_1,j+s_2}}\geq 0$ for all $C_{i+s_1,j+s_2} \in C_{\sum_{i,j}}$;\\[-0.05in]

\textsf{P3}. $m_{k,k}>0$ for all $1\leq k\leq N_1N_2$, equivalently to $C_{i,j}> 0$.\\

\noindent Then we find the condition which guarantees \\[-0.05in]

\textsf{P4}. the non-negativity of the right-hand side $\mathcal{F}$.\\[-0.05in]

At inner points $2\leq i \leq N_1-1$, $2\leq j \leq N_2-1$  from \eqref{EqSD_FD} we get the corresponding coefficients of \eqref{CF} and $\mathcal{F}$
\begin{eqnarray}
&\ds  \hspace{-0.1in}C_{i,j}=\frac{1}{\triangle \tau}+\frac{\sigma_1^2}{h_1^2}+\frac{\sigma_2^2}{h_2^2}-\frac{\sigma_1\sigma_2|\rho'_{i,j}|}{h_1h_2}+r,
\; C_{i\pm1,j}=\frac{\sigma_1^2}{2h_1^2}-\frac{\sigma_1\sigma_2|\rho'_{i,j}|}{2h_1h_2},\;
 C_{i,j\pm 1}=\frac{\sigma_2^2}{2h_2^2}- \frac{\sigma_1\sigma_2|\rho'_{i,j}|}{2h_1h_2},\nonumber\\
&\ds  C_{i-1,j+1}=C_{i+1,j-1}=\frac{\sigma_1\sigma_2\rho'^-_{i,j}}{2h_1h_2}, \ \ \ \  C_{i-1,j-1}=C_{i+1,j+1}=\frac{\sigma_1\sigma_2\rho'^+_{i,j}}{2h_1h_2},\nonumber\\ [-0.1in]
& \label{TP_1} \\ [-0.1in]
&\ds f_{i,j}=\frac{1}{\triangle \tau}u_{i,j}+A_{1}^+\Lambda_{1_{i,j}}^+\frac{u_{i+1,j}-u_{i,j}}{h_1}-A_{1}^-\Lambda_{1_{i,j}}^-\frac{u_{i,j}-u_{i-1,j}}{h_1}\nonumber\\
&\ds +A_{2}^+\Lambda_{2_{i,j}}^+\frac{u_{i,j+1}-u_{i,j}}{h_2}-A_{2}^-\Lambda_{2_{i,j}}^-\frac{u_{i,j}-u_{i,j-1}}{h_2},\nonumber
\end{eqnarray}

Properties \textsf{P1} - \textsf{P3} are fulfilled, owing to \eqref{RH}. We have $|C_{i,j}|- \hspace{-0.3in}\sum\limits_{C_{i+s_1,j+s_2} \in C_{\sum_{i,j}}}\hspace{-0.3in}|C_{i+s_1,j+s_2}|=\frac{1}{\tau}+r$,
 $C_{i,j}\geq \frac{1}{\tau}+r >0$ and all $C_{i+s_1,j+s_2} \in C_{\sum_{i,j}}$ are non-negative. To ensure the property \textsf{P4} we require
 $$
 \frac{1}{\triangle \tau}-\frac{A_{1}^+\Lambda_{1_{i,j}}^+}{h_1}-\frac{A_{1}^-\Lambda_{1_{i,j}}^-}{h_1}
 -\frac{A_{2}^+\Lambda_{2_{i,j}}^+}{h_2}-\frac{A_{2}^-\Lambda_{2_{i,j}}^-}{h_2}\geq 0,
 $$
 which leads to restriction \eqref{RT}.

Let for instance $\partial \Omega'_E \subseteq \partial \Omega'_1$. Thus from \eqref{EqSD_E_FD} we have
 \begin{eqnarray}
&\ds \hspace{-0.1in}C_{N_1,j}=\frac{1}{\triangle \tau}+\frac{\sigma_1^2}{h_1^2}+\frac{\sigma_2^2}{h_2^2}-\frac{\sigma_1\sigma_2|\rho'_{N_1,j}|}{h_1h_2}+r,
\ \ \ \  C_{N_1-1,j}=\frac{\sigma_1^2}{h_1^2}-\frac{\sigma_1\sigma_2|\rho'_{N_1,j}|}{h_1h_2},\nonumber\\
&\ds C_{N_1,j\pm 1}=\frac{\sigma_2^2}{2h_2^2}- \frac{\sigma_1\sigma_2|\rho'_{N_1,j}|}{2h_1h_2},\ \ \ \
 C_{N_1-1,j\pm 1}=\frac{\sigma_1\sigma_2|\rho'_{N_1,j}|}{2h_1h_2}, \nonumber\\ [-0.1in]
& \label{TP_2} \\ [-0.1in]
&\ds f_{N_1,j}=\frac{1}{\triangle \tau}u_{N_1,j}+\left(A_{1}^++\frac{\sigma_1^2}{h_1}-\frac{\sigma_1\sigma_2|\rho'_{N_1,j}|}{h_2}\right)\widehat{g}'_{1_{N_1,j}}+\frac{\sigma_1\sigma_2\rho'^+_{N_1,j}}{h_2}\widehat{g}'_{1_{N_1,j+1}}+\frac{\sigma_1\sigma_2\rho'^-_{N_1,j}}{h_2}\widehat{g}'_{1_{N_1,j-1}}\nonumber\\ & \ds -A_{1}^-\Lambda_{1_{N_1,j}}^-\frac{u_{N_1,j}-u_{N_1-1,j}}{h_1}
 +A_{2}^+\Lambda_{2_{N_1,j}}^+\frac{u_{N_1,j+1}-u_{N_1,j}}{h_2}-A_{2}^-\Lambda_{2_{N_1,j}}^-\frac{u_{N_1,j}-u_{N_1,j-1}}{h_2},\nonumber
\end{eqnarray}
 As before \textsf{P1} - \textsf{P3} follows from \eqref{RH}. The right-hand side is non-negative if additionally to \eqref{RH} we have
 $$
 \frac{1}{\triangle \tau}-\frac{A_{1}^-\Lambda_{1_{N_1,j}}^-}{h_1}
 -\frac{A_{2}^+\Lambda_{2_{N_1,j}}^+}{h_2}-\frac{A_{2}^-\Lambda_{2_{N_1,j}}^-}{h_2}\geq 0 \ \ \hbox{and therefore restriction}\ \ \eqref{RT}.
 $$
 From equations \eqref{EqSD_W_FD}, \eqref{EqSD_N_FD} and \eqref{EqSD_S_FD} we obtain similar results.

 Consider now the corner node $i=N_1$, $j=N_2$, $\{\partial \Omega'_N,\partial \Omega'_E\} \subseteq \partial \Omega'_1$. From \eqref{EqSD_NE_FD} we determine
  \begin{eqnarray}
&\ds \hspace{-0.1in}C_{N_1,N_2}=\frac{1}{\triangle \tau}+\frac{\sigma_1^2}{h_1^2}+\frac{\sigma_2^2}{h_2^2}-\frac{\sigma_1\sigma_2|\rho'_{N_1,N_2}|}{h_1h_2}+r,
\ \ \ \  C_{N_1-1,N_2}=\frac{\sigma_1^2}{h_1^2}-\frac{\sigma_1\sigma_2|\rho'_{N_1,N_2}|}{h_1h_2},\nonumber \\
&\ds C_{N_1,N_2- 1}=\frac{\sigma_2^2}{h_2^2}- \frac{\sigma_1\sigma_2|\rho'_{N_1,N_2}|}{h_1h_2},\ \ \ \
 C_{N_1-1,N_2- 1}=\frac{\sigma_1\sigma_2|\rho'_{N_1,N_2}|}{h_1h_2}, \nonumber 
 \end{eqnarray}
  \begin{eqnarray}
&\ds f_{N_1,N_2}=\frac{1}{\triangle \tau}u_{N_1,N_2}+\left(A_{1}^++A_{2}^++\frac{\sigma_1^2}{h_1}+\frac{\sigma_2^2}{h_2}-\frac{\sigma_1\sigma_2|\rho'_{N_1,j}|}{h_2}-\frac{\sigma_1\sigma_2|\rho'_{N_1,j}|}{h_1}\right)\widehat{g}'_{1_{N_1,N_2}}\nonumber\\ [-0.1in]
& \label{TP_3}\\ [-0.1in]
&\ds +\frac{\sigma_1\sigma_2}{2h_2}\left[\left(|\rho'_{N_1,N_2}|+\rho'^-_{N_1,j}\right)\widehat{g}'_{1_{N_1,N_2-1}}+\rho'^+_{N_1,N_2}\widehat{g}'_{1_{N_1,N_2+1}}\right] \nonumber\\
&\ds +\frac{\sigma_1\sigma_2}{2h_1}\left[\left(|\rho'_{N_1,N_2}|+\rho'^-_{N_1,j}\right)\widehat{g}'_{1_{N_1-1,N_2}} +\rho'^+_{N_1,N_2}\widehat{g}'_{1_{N_1+1,N_2}}\right] \nonumber\\
& \ds -A_{1}^-\Lambda_{1_{N_1,N_2}}^-\frac{u_{N_1,N_2}-u_{N_1-1,N_2}}{h_1}
-A_{2}^-\Lambda_{2_{N_1,N_2}}^-\frac{u_{N_1,N_2}-u_{N_1,N_2-1}}{h_2},\nonumber
\end{eqnarray}
 Evidently, restrictions \eqref{RH} and \eqref{RT} guarantees properties \textsf{P1} - \textsf{P4}. Similar considerations can be applied for \eqref{EqSD_WN_FD}, \eqref{EqSD_SE_FD} and \eqref{EqSD_SW_FD}. \hfill $\Box$

The next results concern the stability of the presented numerical method.

\begin{thm} If $\partial \Omega_1\equiv \emptyset$ or $\partial \Omega_1\not \equiv \emptyset$ and $g_1=0$, $g_s\geq 0$, $s=0,2$ both \eqref{RH} and \eqref{RT} hold,
  then the numerical solution of the problem \eqref{EqT}-\eqref{BC2T} (respectively \eqref{Eq}-\eqref{BC2}), obtained by \eqref{EqSD_FD}, associated with Dirichlet boundary conditions and discretization \eqref{EqSD_W_FD}-\eqref{EqSD_SW_FD} (depending on $\partial\Omega$) is stable (in  maximal discrete norm) with respect to the initial and boundary conditions.
\end{thm}

\noindent\textbf{Proof.} Without loss of generality we will consider \eqref{TP_1}, \eqref{TP_2} and \eqref{TP_3}. The estimates for the other part of the boundary are similar. Let $\|u\|:=\max\limits_{i,j}|u_{i,j}|$. Taking into account restrictions \eqref{RH} and \eqref{RT}, from \eqref{CF} and \eqref{TP_1} we estimate
\begin{equation}\label{IP}
\|\widehat{u}\|\leq \frac{1}{1+r\triangle \tau}\|u\|.
\end{equation}
Similarly, from \eqref{CF}, \eqref{TP_2} and \eqref{TP_3} we again obtain \eqref{IP}.

For homogeneous Neumann boundary conditions we apply the same considerations and after time integration procedure we set
$$
\hspace{1.95in} \|u\|\leq \max\{\|g'_0\|,T\max\limits_{\partial \Omega'_2}g'_2\}.\hspace{1.95in} \Box
$$
\begin{thm} If  $g_s\geq 0$, $s=0,1,2$, $g_1\neq 0$, $\partial \Omega_1\not \equiv \emptyset$,  \eqref{RH}, \eqref{RT} hold  
  then the numerical solution of the problem \eqref{EqT}-\eqref{BC2T} (respectively \eqref{Eq}-\eqref{BC2}), obtained by \eqref{EqSD_FD}, associated with Dirichlet boundary conditions and discretization \eqref{EqSD_W_FD}-\eqref{EqSD_SW_FD} (depending on $\partial\Omega$) is stable (in  maximal discrete norm) with respect to the initial and boundary conditions.
\end{thm}

\noindent\textbf{Proof.} Again we consider \eqref{TP_1}, \eqref{TP_2} and \eqref{TP_3}. As before, at inner points we obtain the estimate \eqref{IP}.  From \eqref{CF}, \eqref{TP_2} and \eqref{TP_3}, substituting $\frac{\sigma_1^2}{h_1}\widehat{g}'_{1_{N_1,j}}=\frac{\sigma_1^2}{h_1}\widehat{u}_{\mathring{x}_{1_{N_1,j}}}$, $\left(\frac{\sigma_1^2}{h_1}+\frac{\sigma_1^2}{h_2}\right)\widehat{g}'_{1_{N_1,N_2}}=\frac{\sigma_1^2}{h_1}\widehat{u}_{\mathring{x}_{1_{N_1,N_2}}}+\frac{\sigma_1^2}{h_2}\widehat{u}_{\mathring{x}_{2_{N_1,N_2}}}$ in view of \eqref{NB} and \eqref{EB}, we get
\begin{eqnarray*}
&&\|\widehat{u}\|\leq \frac{1}{1+r\triangle \tau}\|u\|+\triangle \tau A_{1}^+\|\widehat{g'}_1\|, \\  &&\|\widehat{u}\|\leq \frac{1}{1+r\triangle \tau}\|u\|+\triangle \tau (A_{1}^++A_{2}^+)\|\widehat{g'}_1\|.
\end{eqnarray*}
 Then, taking into account also the Dirichlet boundary conditions (if any), the time integration procedure in general case leads to
$$
\hspace{1.0in} \|u\|\leq \max\{\|g'_0\|+C\max\limits_{\partial \Omega'_1}g'_1, T\max\limits_{\partial \Omega'_2}g'_2\}, \ \ \hbox{where} \ \ C = T(|A_{1}|+|A_{2}|). \ \ \ \ \ \ \ \ \ \ \ \  \Box
$$

\section{Numerical Examples}

 In this section we test the accuracy, convergence rate and  positivity preserving  of the presented numerical methods for model problem \eqref{EqT}-\eqref{BC2T} (and \eqref{Eq}-\eqref{BC2}). Model parameters are $D_1=0.0487902$, $D_2=0$, $\sigma_1=\sigma_2=0.2$, $r=0.0953102$ \cite{Top1}. In agreement with \eqref{RH} we can  choose  $h=h_1=h_2$ ($N=N_1=N_2$).

 When we deal with exact solution (Example 1), the convergence rate in maximal  discrete norm is computed using two consecutive meshes:
 \begin{eqnarray*}
 CR_\infty=\log_2\frac{E^{N/2}_\infty}{E^N_\infty}, \ \ \ \
 E^N_\infty=\max\limits_{{1\leq i,j \leq N_1}}|E_{i,j}^N|,
 \end{eqnarray*}
 where $E_{i,j}^N$ is the difference between the exact and the numerical solutions at point $(x_{1_i},x_{2_j},T)$ on a mesh with $N\times N$ grid nodes in space.

 Alternatively, if the exact solution is not available (Example 2), the  convergence rate is computed by the same formula but now $E_{i,j}^{N}$ is the difference  between two \textit{numerical} solutions, computed on meshes with $N$ and $2N$ grid nodes respectively.

In order to avoid division by zero in uniform flow regions, we add $\varepsilon <<1$ ($\varepsilon=10^{-30}$) to both numerator and
denominator of the gradient ratio \eqref{ratio}.
\vspace{0.1in}

\noindent\textbf{Example 1} (\textsf{Exact solution test}) In the right hand side of the equation \eqref{EqT} we add an appropriate residual function and consider non-homogeneous Neumann boundary conditions on East, North and South boundary ($\partial \Omega'_1\equiv \partial \Omega'_E\cup \partial \Omega'_N \cup \partial \Omega'_S$) and Dirichlet boundary conditions on the West  boundary ($\partial \Omega'_2\equiv \partial \Omega'_W$) such that
\begin{equation*}
u(x_1,x_2,\tau)=e^{-\tau/2}\cos(\pi x_1/3)\cos(\pi x_2/3),
\end{equation*}
is the exact solution of the modified problem \eqref{EqT}-\eqref{BC2T}.
The computations are performed in two domains:
$$
\overline{\Omega}'^A =[-1,1]\times[-1,1],\ \
\ \ \overline{\Omega}'^B \simeq [-\ln(200),\ln(200)]\times[-\ln(200),\ln(200)].
$$
 for $T=0.5$ and  fixed for all time levels  time step $\triangle\tau=h^2$. The results for different values of $\rho_1$, $\rho_2$ in each  domain $\overline{\Omega'}^A$ and $\overline{\Omega'}^B$ are given in Table \ref{t1}.  We observe second-order convergence rate of the numerical method.
\begin{table}[h]
\caption{Errors and convergence rates, Example 1}%
 \label{t1}
\begin{tabular*}{\textwidth}{@{\extracolsep{\fill}}rccccccccc}\hline\\[-0.1in]
&&\multicolumn{2}{c}{$\overline{\Omega}'^A$}&&\multicolumn
{4}{c}{$\hspace{0.6in}\overline{\Omega}'^B$}\\\cline{3-4}\cline{6-10}\\[-0.1in]
\multicolumn {1}{c}{$N$}&&\multicolumn{5}{c}{$\rho_1=-0.2,\; \rho_2=0.6$}&&
\multicolumn{2}{c}{$\rho_1=-1,\; \rho_2=1$}\\\cline{3-7}\cline{9-10}\\[-0.1in]
\multicolumn {1}{c}{}&&\multicolumn {1}{c}{$E^N_\infty$}
 & $CR_\infty$&&\multicolumn {1}{c}{$E^N_\infty$}
 &  $CR_\infty$&&\multicolumn {1}{c}{$E^N_\infty$}
 &  $CR_\infty$
\\[0.03in]
\hline \hline \\[-0.1in]
21&& 6.48015e-4  &    && 1.69489e-2  &  && 1.69709e-2    & \\
41&& 1.58029e-4  & 2.0359   && 4.83743e-3  & 1.8089 && 4.84391e-3    & 1.8088\\
81&& 3.83190e-5  & 2.0441   &&  1.21792e-3   & 1.9898 && 1.21971e-3    & 1.9896\\
161&& 9.38348e-6  & 2.0299   &&  2.86828e-4   & 2.0862 &&  2.87575e-4   &2.0845 \\
321&& 2.32268e-6  & 2.0143   &&  6.84829e-5   & 2.0664 && 6.86652e-5    & 2.0663\\
\hline
\end{tabular*}
\end{table}

\vspace{0.1in}

\noindent\textbf{Example 2} (\textsf{Original problem}) We solve \eqref{EqT}-\eqref{BC2T} (and \eqref{Eq}-\eqref{BC2}) by the  presented numerical method for different initial and boundary conditions.
 All computations are performed in $\overline{\Omega}'^B$ for $\rho_1=-0.2$, $\rho_2=0.6$. For the convergence test we take $\triangle\tau=h^2$  fixed and $T=2$, while the given plots are for
 different time and time steps, satisfying  equality in \eqref{RT}. 
 We denote by $E$ the exercise price,   $w_i$ is  the weight of the  $i$-th asset,  'cap' parameter is used for capped-style options, BS (Price, Strike, Time) is the Black-Scholes vanilla Put/Call option price.

 We consider the following test problems:
\begin{itemize}
\item[\textsf{TP1:}] \textit{European exchange option with pay-off:} $P(S_1,S_2)=\max\{0,S_2-S_1\}$. We use the pay-off function as the source for the Dirichlet condition \cite{KN}. Namely, $\partial \Omega'_1\equiv \emptyset$ and $g_2(S_1,S_2,t)=P(S_1,S_2)$.

\item[\textsf{TP2:}] \textit{Worst-off two Call option with barrier} \cite{ZVF}. Now $P(S_1,S_2)=\max\{0,\min\{S_1,S_2\}-E\}$ and  $\partial \Omega'_1\equiv \emptyset$, $g_2(S_1,S_2,t)=P(S_1,S_2)$.

\item[\textsf{TP3:}] \textit{Capped  Put on a basket of two equities} \cite{Top0,Top1}. The initial function is $g_0=\min\{\textnormal{cap},\max\{0,E-w_1S_1-w_2S_2\}\}$, boundary conditions are \eqref{BC2} ($\partial \Omega'_1\equiv \emptyset$) with
    $$g_2=\left\{\begin{array}{ll}
    0 & \hbox{on} \ \ \partial\Omega_N\cup\partial\Omega_E,\\
    BS(S_1,\frac{E}{w_1},t)-BS(S_1,\textnormal{cap},t)& \hbox{on} \ \ \partial\Omega_S,\\[0.035in]
    BS(S_2,\frac{E}{w_2},t)-BS(S_2,\textnormal{cap},t)& \hbox{on} \ \ \partial\Omega_W,
    \end{array}\right.
    $$
The boundary conditions at $\partial \Omega_W$ and $\partial \Omega_S$ represents the prices of capped European option with strike prices of $E/w_1$ and $E/w_2$, respectively \cite{Top1}.


  \item[\textsf{TP4:}] \textit{Two-asset barrier options} \cite{Hn,Top1}. We consider $\partial \Omega'_1\equiv \partial \Omega'_E\cup \partial \Omega'_N \cup \partial \Omega'_S$,  $\partial \Omega'_2\equiv \partial \Omega'_W$, $g_0=\max\{0, w_1S_1-E\}$, $g_2=0$, $g_1=0$ on $\Omega_S\cup  \Omega_N$, $g_1=1$ on $\Omega_E$.

       \item[\textsf{TP5:}] \textit{Capped  Call on a Basket of two equities} \cite{Top0,Top1}. In this case $\partial \Omega'_1\equiv \partial \Omega'_E\cup \partial \Omega'_N $,  $\partial \Omega'_2\equiv \partial \Omega'_W\cup \partial \Omega'_S$, $g_0=\min\{\textnormal{cap},\max\{0,w_1S_1+w_2S_2-E\}\}$, $g_1=0$ on $\Omega_N\cup  \Omega_E$ and
            $$g_2=\left\{\begin{array}{ll}
        BS(S_1,\textnormal{cap},t)-BS(S_1,\frac{E}{w_1},t)& \hbox{on} \ \ \partial\Omega_S,\\[0.035in]
    BS(S_2,\textnormal{cap},t)-BS(S_2,\frac{E}{w_2},t)& \hbox{on} \ \ \partial\Omega_W,
    \end{array}\right.
    $$
\end{itemize}

In Table \ref{t2} we give convergence rate ($CR_\infty$), computed on three consecutive meshes, for each test problem, $E=100$, $w_1=w_2=1$, cap $=10$.
\begin{table}[h]
\caption{Convergence rates for different problems, $\triangle\tau=h^2$, $T=2$, Example 2}%
 \label{t2}
\begin{tabular*}{\textwidth}{@{\extracolsep{\fill}}ccccccc}\hline\\[-0.1in]
\multicolumn {1}{c}{\textsf{space meshes}}&&\multicolumn {1}{c}{\textsf{TP1}}
 &\multicolumn {1}{c}{\textsf{TP2}}
&\multicolumn {1}{c}{\textsf{TP3}}
 &\multicolumn {1}{c}{\textsf{TP4}} &\multicolumn {1}{c}{\textsf{TP5}}
 \\[0.03in]
\hline \hline \\[-0.1in]
21-41-81 && 1.4458   & 1.3809   &  0.7447   &  1.1625 & 0.7443  \\
41-81-161 && 1.8038  & 1.5757  & 1.4963   &  1.4525  &  1.4732 \\
81-161-321 && 2.0477 & 1.7639  & 1.8234 & 1.8884 &  1.8022\\
\hline
\end{tabular*}
\end{table}
We observe that the order of convergence very close to 2 for all problems \textsf{TP1-TP5}.

\section*{Conclusions} In this paper we develop   second-order in space implicit-explicit finite difference method, based on the van Leer flux-limiter technique, for the worst-case pricing model in financial mathematics.
  Under mild time and space step restrictions the proposed method is stable (with respect to initial and boundary conditions) and preserves the non-negativity of the numerical solution. Van Leer's flux limiter technique is implemented appropriately also for non-homogeneous Neumann boundary conditions, ensuring second order convergence rate and possibility to guarantee the  positivity preserving property of the numerical solution.

  Various numerical examples confirm the theoretical statements and illustrate the second order convergence in space variable.

  The very important question - to find interface curve (in the one dimensional case) or surface (in the two-dimensional case) where the sign of $\Gamma_{cross}$  changes  and on this base to construct numerical method for the corresponding linear problems on both sides of the interface will be the main subject of our next work.

\section*{Acknowledgement}

 This research was supported by the European Union under Grant Agreement number 304617
(FP7 Marie Curie Action Project Multi-ITN STRIKE - Novel Methods in Computational Finance)
and  Bulgarian National Fund of Science under Project   DID 02/37-2009.



\end{document}